\documentclass[a4paper, 12pt]{article}

\usepackage[english]{babel}

\usepackage[a4paper,top=2cm,bottom=3cm,left=3cm,right=3cm,marginparwidth=1.75cm]{geometry}
 \usepackage[utf8]{inputenc}
\usepackage[colorlinks=true, allcolors=blue]{hyperref}

\usepackage{amsmath,amsthm,amsfonts,amssymb,amscd}
\usepackage{stmaryrd}
\usepackage{enumitem}
\usepackage{mathrsfs}
\usepackage{longtable,geometry}
\usepackage{mathtools}
\usepackage{yhmath}
\usepackage{pgf,tikz,pgfplots}
\pgfplotsset{compat=1.15}
\usetikzlibrary{arrows}
\usepackage{dsfont}
\usepackage{xcolor}
\usepackage{graphicx}
\usepackage{listings}
\usepackage{hyperref}
\usepackage{stmaryrd}
\usepackage{verbatim}
\usepackage{subcaption}
\usepackage{float}
\usepackage{todonotes}
\usepackage{nicematrix}
\usepackage{lineno}

\usepackage{appendix}

\newcommand{\R}{\mathbb{R}}

\newcommand{\N}{\mathbb{N}}

\newcommand{\Haus}{\mathcal{H}^{N-1}}

\newcommand{\mres}{\mathbin{\vrule height 1.6ex depth 0pt width 0.13ex\vrule height 0.13ex depth 0pt width 1.3ex}}

\title{Optimal regularity up to the boundary for Plateau-quasi-minimizers}
\author{Eve Machefert \footnote{INSA Lyon, CNRS, Ecole Centrale de Lyon, Université Claude Bernard Lyon 1, Université Jean
Monnet, ICJ UMR 5208, 69621 VILLEURBANNE, FRANCE, eve.machefert@insa-lyon.fr}}

\newtheorem{thm}{Theorem}[section] 

\newtheorem{definition}[thm]{Definition} 
\newtheorem{prop}[thm]{Proposition}

\newtheorem{rmq}[thm]{Remark}
\newtheorem{ex}[thm]{Example}
\newtheorem{lemme}[thm]{Lemma}

\DeclareMathOperator\spt{spt}

\numberwithin{equation}{section}

\begin{document}


\maketitle

\begin{abstract}
We study the regularity of quasi-minimal sets (in the sense of David and Semmes) with a boundary condition, which can be interpreted as quasi-minimizers of Plateau's problem in co-dimension one. For these Plateau-quasi-minimizers, we establish the optimal regularity, which is a characterization by bi-John domains with Ahlfors regular boundary. This requires to investigate the Ahlfors regularity and also the uniform rectifiability of those sets, up to the boundary.    
\end{abstract}

\tableofcontents

\section{Introduction}
 
The study of minimal sets is a major topic in geometric measure theory and variational problems involving an interface. In particular, the notorious Plateau problem aims to find a surface that minimizes its area and spans a given curve. One application is to determine the shape of soap films that span a prescribed contour. In the case of co-dimension one, Plateau's problem can be solved by De Giorgi's approach, using the theory of sets of finite perimeter (see \cite{de2017direct}, \cite{DeGiorgi1954}, \cite{DeGiorgi1955}, \cite{DeGiorgi1961Complementi}, \cite{giusti1984minimal} and \cite{maggi_sets_2012}). This approach consists in minimizing among sets of finite perimeter that satisfy a boundary condition instead of surfaces that spans a curve. In this case, the optimal surface is represented by the essential boundary $ \partial^* \Omega_0$ of the minimal set $\Omega_0$. 

This paper focuses on the study of quasi-minimizers of perimeter, a wider class than the one of perimeter minimizers. The concept of quasi-minimizers of perimeter was introduced in \cite{david1998quasiminimal} by David and Semmes. They proved that theses quasi-minimizers are characterized by bi-John domain with regular boundary (in a sense of Definition~\ref{def : regbdd}). Later in \cite{labourie2024optimalregularityquasiminimalsets}, this regularity was extended to quasi-minimal sets without the assumption that its boundary separates the space in only two connected components.

As these previous works do not consider a boundary condition in their study, this paper generalizes the regularity result of \cite{david1998quasiminimal} by adding such condition. This generalization establishes the optimal regularity up to the boundary of perimeter quasi-minimizers with a boundary condition, which can be interpreted as quasi-minimizers of Plateau's problem in co-dimension one. Specifically, we characterize quasi-minimizers that satisfy a boundary condition by bi-John domain with Ahlfors regular boundary.

John domains were introduced in 1961 by Fritz John (see \cite{john1961rotation}) in connection with his work on elasticity, and studied by Martio–Sarvas (see\cite{martio1979injectivity}). Bi-John domains are defined as John domains such that their complementary set is also a John domain.

\medskip

\textbf{Plateau-quasi-minimizer and main result.}
In this paper we introduce the notion of Plateau-quasi-minimizers, based on quasi-minimizers of perimeter, introduced by David and Semmes in \cite{david1998quasiminimal} and studied in patricular in \cite{rigot2000uniform}, to which we add a boundary constraint, denoted by $E_0$. This added constraint can be interpreted as the boundary condition of Plateau's problem in the De Giorgi formulation, see \cite[(12.29)]{maggi_sets_2012}. We work in a open bounded domain, denoted by $\mathcal{D}$, with $\mathcal{C}$ an open convex non empty set in $\mathcal{D}$. We start by defining the class of competitors.

 \begin{definition}[Competitor]
\label{competitors}
A Borel set $\Omega \subset \mathcal{D}$ is a competitor for Plateau's problem \eqref{def : pblimit} if it is of finite perimeter in $\mathcal{D}$ (\emph{i.e.}, $\chi_{\Omega} \in BV(\mathcal{D})$) and satisfies the boundary condition : $\Omega \setminus \overline{\mathcal{C}} = E_0 \setminus \overline{\mathcal{C}}$.
\end{definition}
The Plateau problem associated with the boundary condition $E_0$ is defined as follows.

\begin{definition}[Plateau's problem]
 \begin{equation}
 \label{def : pblimit}
     \min \{ P(\Omega ,\mathcal{D}) , \Omega \text{ a competitor as defined in Definition~\ref{competitors} } \}.
 \end{equation} 
\end{definition}
This problem admits solutions, as recalled in Section~\ref{Section existence}. Now, let us introduce the main object of our study: Plateau-quasi-minimizers.

\begin{definition}[Plateau-quasi-minimizer]
\label{def quasi min}
A competitor of \eqref{def : pblimit} (see Definition~\ref{competitors}) $\Omega_{0}$ is called a Plateau-quasi-minimizer (or P-quasi-minimizer for short) if there exists a constant $Q\geqslant 1$ such that, for all $\Omega$ competitor of \eqref{def : pblimit}, we have
\[\Haus((\partial^*\Omega_0 \setminus \partial^*\Omega)\cap \mathcal{D}) \leqslant Q \Haus((\partial^*\Omega \setminus \partial^*\Omega_0)\cap \mathcal{D}).\]
\end{definition}

Notice that a minimizer of \eqref{def : pblimit} is a Plateau-quasi-minimizer if and only if $Q=1$. However, if $Q>1$, the notion of Plateau-quasi-minimizer is strictly weaker than the notion of minimizer. Additionally, choosing $Q<1$ would lead to a contradiction since we can consider a $\Omega = \Omega_0$.

In the following, we will also refer to the essential boundary of a Plateau-quasi-minimizer as a quasi-optimal-surface.

\begin{rmq}
\label{rq : quasimin}
Let $\Omega_0$ be a Plateau-quasi-minimizer and $\Omega$ a competitor. If $x \in \partial^*\Omega_0$ and $r>0$ are such that $\Omega_0 \Delta \Omega \subset B(x,r)$, then Definition~\ref{def quasi min} yields 
\[P(\Omega_0, B(x,r)) \leqslant Q P(\Omega, B(x,r)).\]
\end{rmq}

To state our Theorem, we recall first the definition of John domains. 

\begin{definition}[John domain]
\label{def John}
Let $\Omega \subset \mathcal{D}$ be an open bounded set. $\Omega$ is a John domain in $\mathcal{D}$, if there exist a point $z_0 \in \Omega$, called the center, and a constant $C>0$ such that for all $x \in \Omega$ there is a Lipschitz curve $\alpha : [0,|x-z_0|] \rightarrow \Omega$ connecting $x$ to the center $z_0$ ($\alpha(0) = x$ and $\alpha(|x-z_0|)=z_0$) and 
\[d(\alpha(t),\R^{N} \setminus \Omega) \geqslant C^{-1}t, \text{ for } 0 \leqslant t \leqslant |x-z_0|.\]
Moreover, we say that $\Omega$ is a bi-John domain in $\mathcal{D}$ if both $\Omega$ and $\mathcal{D} \setminus \overline{\Omega}$ are John domains.
\end{definition}

Roughly, a John domain is an open set, path-connected, such that any path connecting a point to the center moves away from the boundary at a rate that is at least linear. Therefore, John domains cannot have cusps, for instance.

The goal of this paper is to establish the optimal regularity for these Plateau-quasi-minimizers, under some assumptions on the boundary constraint discussed below. Here is our main result. 

\begin{thm}[Characterization of Plateau-quasi-minimizer by bi-John domains]
\label{th john dom intro} Let $(\mathcal{D}, \mathcal{C}, E_0)$ satsifying Hypothesis H (see Definition~\ref{def H}) and $\Omega_0$ be a competitor, as in Definition~\ref{competitors}.
Then, $\Omega_0$ is a Plateau-quasi-minimizer if and only if there exists an equivalent open set $\Omega$, which is a bi-John domain with Ahlfors-regular boundary, such that $\spt \mu_\Omega =\partial \Omega$.  
\end{thm}

In other words, our main result states that if we replace a quasi-minimizer $E_0$ by another one in a convex subset $\mathcal{C}$, the resulting set, by gluing the two, remains a quasi-minimizer. Of course, this is not true in general but it is valid under some natural conditions, such as Condition B' below, on how $E_0$ meets $\mathcal{C}$.

As mentionned before we extend to the boundary the interior regularity known for quasi-minimizer. It has been proven in \cite{david1998quasiminimal} that quasi-minimizers (without boundary constraint) are bi-John domain with Ahlfors regular boundary and satisfy Condition B (see Definition~\ref{condition B}), introduced in \cite[Definition 2, p.102]{david1988morceaux} and studied in \cite{rigot2000uniform}. Naturally, we cannot hope that this regularity result holds without assuming that the boundary constraint $E_0$ itself satisfies this regularity outside $\mathcal{C}$, which is insured by Hypothesis H given below. For sake of clarity, we introduce the notation $\Sigma := \partial E_0\cap (\mathcal{D}\setminus \overline{\mathcal{C}}) $. Let us define the assumptions on the boundary constraint $E_0$.

\begin{enumerate}[label = (H\arabic*)]
\item
\label{H2} 
 $\Sigma$ is {\bf Ahlfors regular}, \emph{i.e.}, there exist two constants $C_1,C_2>0$ such that for all $x\in \Sigma$ and all $r>0$ such that 
$B(x,r) \subset \mathcal{D}$, 
\begin{equation*}
C_1 r^{N-1} \leqslant \Haus(\Sigma \cap B(x,r))  \leqslant C_2 r^{N-1}.
\end{equation*}
\item
\label{H3}$E_0$ satisfies {\bf Condition B'}, defined in Definition~\ref{condition B'}. 
\item
\label{H1} $E_0 \setminus \overline{\mathcal{C}}$ and $(\mathcal{D} \setminus \overline{E_0}) \setminus \overline{\mathcal{C}}$ are {\bf domains of isoperimetry} in $\mathcal{D}$, as defined in Definition~\ref{def : domain isop}.

\end{enumerate}

Condition B', defined as follows, is a Condition B, as introduced in \cite[Definition 2, p.102]{david1988morceaux}, relative to the convex set $\mathcal{C}$. 

\begin{definition}[Condition B']
\label{condition B'} Let $E \subset \mathcal{D}$ be an Borel set of finite perimeter in $\mathcal{D}$. $E$ satisfies Condition B' if it is open, $\partial E$ is Ahlfors regular in $\mathcal{D}$ and there exists a constant $C > 1$ such that, for all $x \in  \partial E \setminus \overline{\mathcal{C}}$ and $r>0$ satisfying $B(x,r) \subset \mathcal{D}$, there exist two balls $B_1 \subset (B(x,r) \cap E)\setminus \overline{\mathcal{C}}$ and $B_2 \subset (B(x,r) \setminus \overline{E})\setminus \overline{\mathcal{C}}$ with radii larger or equal to $C^{-1}r$. 
\end{definition} 

We can now define Hypothesis H, needed for the main Theorem of this paper. 

\begin{definition}[Hypothesis H]
\label{def H}
We say that $(\mathcal{D}, \mathcal{C}, E_0)$ satsifies Hypothesis H if $\mathcal{D}$ is an open bounded domain, $\mathcal{C}$ is an open convex non empty set in $\mathcal{D}$ and $E_0$ is a Borel set of finite perimeter in $\mathcal{D}$ satisfying Hypotheses~\ref{H2},\ref{H3} and \ref{H1}, and such that $\partial^* E_0 = \partial E_0$ (so in particular $\spt \mu_{E_0} = \partial E_0$).
\end{definition}

\textbf{Necessity of Hypothesis H.} As mentionned above, the assumptions on $E_0$ are natural considering the regularity we want to establish. Notice that John domains are known to be domains of isoperimetry. Therefore Hypothesis~\ref{H1} is weaker than assuming that outside $\mathcal{C}$, $E_0$ is a bi-John domain. This assumption is however not unexpected since, in \cite{david1998quasiminimal}, it is proven that quasi-minimizers of perimeter are domains of isoperimetry.  

\begin{figure}[h!]
  \centering
    \subfloat[\centering Boundary constraint not satisfying the hypothese~\ref{H3} and \ref{H1}]{\includegraphics[scale = 0.34]{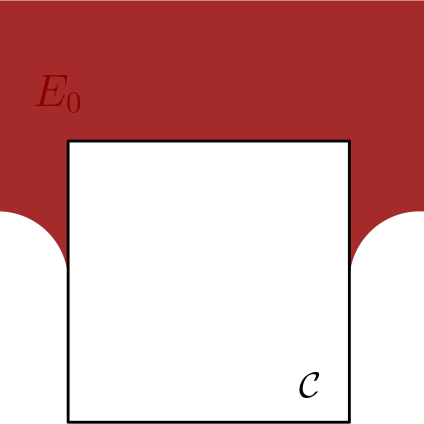} \label{fig:contre-ex1E0}} \quad \quad 
    \subfloat[\centering Plateau-quasi-minimizer, which is not a John domain]{\includegraphics[scale = 0.34]{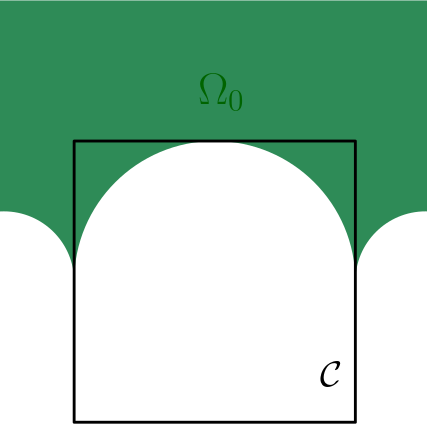} \label{fig:contre-ex1O0}}
  \caption{ Example to justify hypotheses~\ref{H3} and \ref{H1}}
  \label{fig:contre-ex1}
\end{figure}

\begin{figure}
  \centering
    \subfloat[\centering Alternative boundary constraint ]{\includegraphics[scale = 0.34]{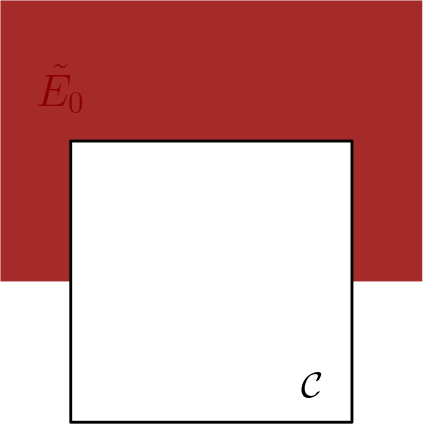} \label{fig:contre-ex2E0}} \quad \quad 
    \subfloat[\centering Plateau-quasi-minimizer for the alternative boundary constraint]{\includegraphics[scale = 0.34]{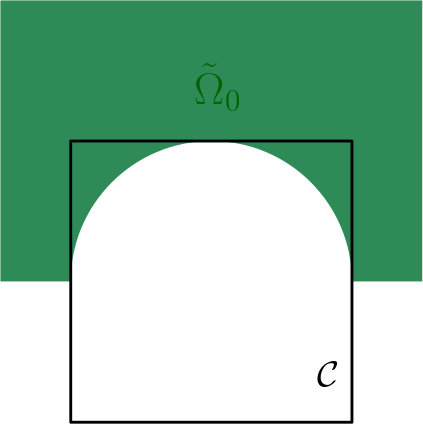} \label{fig:contre-ex2O0}}
  \caption{ Example to justify that $\Omega_0$ is a Plateau-quasi-minimizer}
  \label{fig:contre-ex2}
\end{figure}

In Section~\ref{section cylindre}, we present an example of a set $E_0$ satisfying these assumptions. Let us here introduce examples to justify the necessity of these assumptions on the boundary condition. Consider $\mathcal{C}$ to be a square in 2D and the boundary constraint $E_0$ in red in Figure~\ref{fig:contre-ex1E0}.
 Such $E_0$ does not satisfy Condition B' because it forms a cusp outside $\mathcal{C}$, and for the same reason it is not a John domain and thus not a domain of isoperimetry. Hence, Hypotheses~\ref{H1} and ~\ref{H3} are not verified by this $E_0$. We will now justify why $\Omega_0$ defined in Figure~\ref{fig:contre-ex1O0} is a Plateau-quasi-minimizer. Consider a different boundary condition $\tilde{E_0}$ as represented in Figure~\ref{fig:contre-ex2E0}. Then $\tilde{\Omega}_0$, in green in Figure~\ref{fig:contre-ex2O0}, is clearly a bi-John domain with regular boundary. Hence from the result of David and Semmes, it is a quasi-minimizer inside $\mathcal{C}$. This yields that $\Omega_0$ is a quasi-minimizer inside of $\mathcal{C}$, and since outside of $\overline{\mathcal{C}}$ all competitors are equal, there only remains to verify the condition at the boundary $\partial \mathcal{C}$. However, from the definition of $\Omega_0$ and $\tilde{\Omega}_0 $ it is clear that $\Haus(\partial ^* \Omega_0 \cap \partial \mathcal{C})=0$. Thus, $\Omega_0$ is a Plateau-quasi-minimizer.
 Again, due to the cusps formed at the boundary of the square, $\Omega_0$ is not a John domain.

Notice also that for any set of finite perimeter $E$, there exists a Borel set of finite perimeter $F$ equivalent to $E$ such that $\spt \mu_F = \partial F$ (see \cite[Proposition 12.19]{maggi_sets_2012}). Therefore the assumption that $\spt \mu_{E_0} = \partial E_0$ is not restrictive. 

\medskip

As briefely mentionned in Theorem~\ref{th john dom intro}, the regularity is proved on a equivalent set.  
\begin{rmq}
\label{rq : equivalent sets}
Recall that two sets $E$ and $F$ are equivalent if $| E\Delta F|=0$. Therefore, two equivalent sets have the same density points, and thus the same essential boundary. As recalled above, for any set of finite perimeter $E$, there exists a Borel set of finite perimeter $F$ equivalent to $E$ such that $\spt \mu_F = \partial F$ (see \cite[Proposition 12.19]{maggi_sets_2012}). In particular, this implies that $\overline{\partial_*F}=\partial F$ (see \cite[Remark 15.3]{maggi_sets_2012}), and, since $\partial_*F \subset \partial^*F \subset \partial F$, it also yields that $\overline{\partial^*F}=\partial F$. Let us justify that, if $E$ is a competitor, \emph{i.e.}, it satisfies the boundary condition then, $F$ is also competitor. Indeed, by definition $F$ coincides almost everywhere with $E_0$ in $\mathcal{D}\setminus \overline{\mathcal{C}}$. Moreover, both $E_0$ and $F$ have a topological boundary that coincides with the closure of their essential boundary. This leads to $F$ satsisfying the boundary condition. Hence, any set equivalent to a P-quasi-minimizer (respectively minimizer) of \eqref{def : pblimit} such that the support of the Gauss-Green measure coincides with the closure of its essential boundary is itself a P-quasi-minimizer (respectively minimizer) of \eqref{def : pblimit}. In the rest of the paper, we will thus prove regularity results for P-quasi-minimizers by establishing the desired regularity for an equivalent set satisfying this condition.
\end{rmq}

\textbf{Comments on the proof.} To make such comments, we must first provide some additional definitions.
\begin{definition}[Ahlfors regularity in an open set]
\label{def AR}
Let $U$ be an open set and $S \subset U$ a closed subset. $S$ is said to be Ahlfors regular in $U$ if there exist two constants $C_1>0$ and $C_2>0$ such that, for all $x\in S \cap U$ and all $r>0$ satisfying $B(x,r) \subset U$, 
\[C_1 r^{N-1} \leqslant \Haus(S\cap B(x,r)) \leqslant C_2 r^{N-1}.\]
\end{definition}

This allows us to define the notion of regular boundary.

\begin{definition}[Regular boundary]
\label{def : regbdd}
Let $\Omega_0$ be a Borel set of finite perimeter in $\mathcal{D}$. We shall say that its boundary is regular if $\Omega_0$ is open, $\spt \mu_{ \Omega_0} = \partial \Omega_0$, and $\partial \Omega_0$ is Ahlfors regular in $\mathcal{D}$.

\end{definition}

\begin{rmq}
As explained before, for any set of finite perimeter $E$, there exists a Borel set of finite perimeter $F$ equivalent to $E$ such that $\spt \mu_{ F} =\partial F$. Therefore, Definition~\ref{def : regbdd} assumes that this equivalent set is open with Ahlfors regular topological boundary.

\end{rmq}

To prove Theorem~\ref{th john dom intro}, we will use the fact that Plateau-quasi-minimizers satisfy the so called Condition B (see Definition~\ref{condition B}) introduced in \cite[Definition 2, p.102]{david1988morceaux}.

\begin{definition}[Condition B]
\label{condition B} Let $E \subset \mathcal{D}$ be a Borel set of finite perimeter in $\mathcal{D}$. $E$ satisfies Condition B if it is open, $\partial E$ is Ahlfors regular in $\mathcal{D}$ and there exists a constant $C > 1$ such that, for all $x \in  \partial E \cap \mathcal{D}$ and $r>0$ satisfying $B(x,r) \subset \mathcal{D}$, there exist two balls $B_1 \subset B(x,r) \cap E$ and $B_2 \subset B(x,r) \setminus \overline{E}$ with radii larger or equal to $C^{-1}r$. 
\end{definition} 

Note that Condition B requires the boundary to be Ahlfors regular. Therefore, the first step of the proof, developped in Section~\ref{section AR}, is to establish the Ahlfors regularity up to the boundary of quasi-optimal surfaces. In other words, we show that there exist two constants $C_{1},C_{2} > 0$ such that, for all $x \in \partial \Omega_{0}\cap \mathcal{D}$ and all $r>0$ satisfying $B(x,r) \subset \mathcal{D}$, we have
\begin{equation}
        C_{1} \leqslant \frac{\Haus(\partial^{*}\Omega_{0}\cap B(x,r))}{r^{N-1}} \leqslant C_{2}.
        \end{equation} 
This formulation can also be rewritten as 
    \begin{equation}
    \label{eq : ARess}
        C_{1} \leqslant \frac{P(\Omega_{0},B(x,r))}{r^{N-1}} \leqslant C_{2}.
        \end{equation}
        
Here, $\partial^*\Omega_0$ is not necessarily closed as required in Definition~\ref{def AR}. However, we will see in the proof of the Ahlfors regularity, that for a Plateau-quasi-minimizer  $\Omega_0$ such that $\spt \mu_{\Omega_0} = \partial \Omega_0$, we have $\Haus$-a-e. $\partial \Omega_0 = \partial^* \Omega_0$. Thus \eqref{eq : ARess} will yield the Ahlfors regularity, in the sense of Definition~\ref{def AR}, for the closed set $\partial \Omega$.

This regularity result is interesting in itself. It is, for instance, used in \cite{BBLM2025} to prove a $\Gamma$-convergence type result, where we present a phase field approximation for Plateau's problem \eqref{def : pblimit}. Thereby, the authors generalize the approximation for Steiner's problem, introduced in \cite{lemenant2014modica,bonnivard2015approximation}. This further motivates the study of Plateau's problem \eqref{def : pblimit}, and of the Ahlfors regularity of its solutions.
        
Then, in Section~\ref{section unif rect}, we prove that Condition B (see Definition~\ref{condition B}) is satisfied by Plateau-quasi-minimizers. In particular, this yields the uniform rectifiability up to the boundary of a quasi-optimal surfaces. And, more precisely, it implies a slightly stronger result: quasi-optimal surfaces have Big Pieces of Lipschitz Graph (BPLG) up to the boundary.

Finally, in Section~\ref{section John}, we establish the optimal regularity result for Plateau-quasi-minimizers. In Theorem~\ref{th john dom intro} we establish the characterization of Plateau-quasi-minimizers by bi-John domain satisfying the boundary condition and with regular boundary, in the sense of Definition~\ref{def : regbdd}.

Lastly, studying the regularity of Plateau-quasi-minimizer instead of solutions of Plateau's problem allows us to establish this regularity not only in a convex set $\mathcal{C}$ but also in bi-Lipschitz images of $\mathcal{C}$. This extension is carried in Section~\ref{section BiLip}.

\medskip

Let us enforce the differences with \cite{david1998quasiminimal} and the main contributions of this article. Notice that in this article, the results are formulated with essential boundaries and perimeters, when David and Semmes used the $\Haus$ measure of topological boundaries. This allowed us to use geometric measure theory results to simplify some computations made in \cite{david1998quasiminimal}. 

Then, as mentionned before, we added a boundary constraint and extended the regularity up to the boundary. 
Firstly, Lemma~\ref{lem AR} is the crucial result that will allow us to extend the existing Ahlfors regularity  and uniform rectifiability result to the boundary for our problem with boundary condition, see Proposition~\ref{th : AhlReg bdd} and Proposition~\ref{prop : rect bdd}. 
Secondly, Lemma~\ref{lem : domain isop} is the key argument to show that a quasi-minimiser is a bi-John domain, see Proposition~\ref{prop : dom isop}. Lastly, for the converse implication, namely to show that a bi-John domain with regular boundary, in the sense of Definition~\ref{def : regbdd} and satisfying the boundary condition is a quasi-minimizer, we defined differently a weighted problem. The main argument for this part is proved in Lemma~\ref{lem John reciproque}. 

In addition, we also develop the full proof of results, presented in \cite{maggi_sets_2012} as two exercices, in Lemma~\ref{exo : 4.1} and Lemma~\ref{Maggi_cvx}.

\medskip

\textbf{Notations.}
In the following, for a Lebesgue-measurable set $E \subset \R^{N}$ we denote by $|E|$ its Lebesgue measure. In particular, we will denote by $\omega_{N} = |B(0,1)|$ the Lebesgue measure of the unit ball in $\R^{N}$. The $k$-Hausdorff measure will be denoted by $\mathcal{H}^{k}$. And when $\mu$ is a Borel measure and $E \subset \R^{N}$ is a Borel set, we denote by $\mu \mres E$ the measure defined as $\mu \mres E(F) = \mu(E \cap F)$. We also denote by $\mathrm{spt}(\mu)$ the support of the measure $\mu$. 

The indicator function of the set $\Omega$ will be denoted by $\chi_{\Omega}$. For $f$ a Lipschitz function, we denote by $Lip(f)$ its Lipschitz constant. And, $\amalg$ represents the disjoint union between sets.

$\partial E$, $\partial^{*}E$ and $\partial_{*}E$ are respectively used to denote the topological boundary, the essential boundary and the reduced boundary of $E$. For $E$ a set of locally finite perimeter the Radon measure $\mu_{E}$ denote the weak derivative of $\chi_{E}$, such that
\[\int_{E} \mathrm{div } (T) = \int_{\R^{N}}T \cdot d\mu_{E}, \hspace{1cm} \forall T \in C_{c}^{1}(\R^{N},\R^{N}).\]
We recall the definition of the measure-theoretic outer unit normal to $E$ denoted by
$\nu_E ~:~ \partial_* E \to \mathbb{S}^{N-1}$ :   
\[\nu_{E}(x) := \lim_{r\to 0^+}\frac{\mu_E(B(x,r))}{|\mu_E|(B(x,r))}, \text{ for } x \in \partial_*E.\]
For better readability, we denote $\mathcal{A}^+ := E_0 \setminus \overline{\mathcal{C}} $ and 
 $\mathcal{A}^- := (\mathcal{D} \setminus \overline{E_0}) \setminus \overline{\mathcal{C}} $. 
\medskip

\textbf{Acknowledgments.}
 I sincerely acknowledge Camille Labourie, as well as Guy David, for the fruitful discussions on this topic. I would also like to thank Antoine Lemenant for introducing this problem to me and for his advices and suggestions regarding this work. Finally, I am grateful to Matthieu Bonnivard for his help with this project. This work was partially supported by the IUF grant of Antoine Lemenant and by the ANR project STOIQUES.

\section{Preliminary}

\subsection{Some classical properties for sets of finite perimeter}

We remind the reader of some results on sets of finite perimeter, which will be useful in the following sections.

\begin{lemme}{\cite[(3.62)]{ambrosio2000oxford}}
For $E$ a set of finite perimeter in the open set $\Omega \subset \R^{N}$, 
\[P(E,\Omega) = \mathcal{H}^{N-1}(\partial ^{*}E \cap \Omega).\]
\end{lemme}

\begin{prop}{\cite[Theorem 16.3]{maggi_sets_2012}}
\label{prop Maggi}
If $E$ and $F$ are sets of locally finite perimeter and we let 
\[\{\nu_{E} = \nu_{F}\} = \{x\in \partial_{*}E\cap \partial_{*}F : \nu_{E}(x) = \nu_{F}(x)\},\]
\[\{\nu_{E} = -\nu_{F}\} = \{x\in \partial_{*}E\cap \partial_{*}F : \nu_{E}(x) = -\nu_{F}(x)\},\]
then $E\cap F$, $E\setminus F$ and $E\cup F$ are sets of locally finite perimeter, with
\begin{align}
&\mu_{E\cap F} = \mu_{E} \mres F^{(1)} + \mu_{F} \mres E^{(1)} + \nu_{E}\mathcal{H}^{N-1} \mres \{\nu_{E} = \nu_{F}\}, \\
&\mu_{E\setminus F} = \mu_{E} \mres F^{(0)} - \mu_{F} \mres E^{(1)} + \nu_{E}\mathcal{H}^{N-1} \mres \{\nu_{E} = -\nu_{F}\}, \\
&\mu_{E\cup F} = \mu_{E} \mres F^{(0)} + \mu_{F} \mres E^{(0)} + \nu_{E}\mathcal{H}^{N-1} \mres \{\nu_{E} = \nu_{F}\}. 
\end{align}
Moreover, for every Borel set $G \subset \R^{N}$,
\begin{align}
&P(E\cap F, G) = P(E,F^{(1)}\cap G) +P(F,E^{(1)}\cap G) + \mathcal{H}^{N-1}(\{\nu_{E} = \nu_{F} \}\cap G), \label{Maggi inter}\\
&P(E\setminus F, G) = P(E,F^{(0)}\cap G) +P(F,E^{(1)}\cap G) + \mathcal{H}^{N-1}(\{\nu_{E} = -\nu_{F} \}\cap G), \label{Maggi prive}\\
&P(E\cup F, G) = P(E,F^{(0)}\cap G) +P(F,E^{(0)}\cap G) + \mathcal{H}^{N-1}(\{\nu_{E} = \nu_{F} \}\cap G).\label{Maggi union}
\end{align}
\end{prop}

The following Lemmas state that intersecting a set of finite perimeter with a convex set reduces its perimeter. This result is presented in  \cite[Exercise 15.13]{maggi_sets_2012} and \cite[Exercise 15.14]{maggi_sets_2012}. For the reader's convenience, we give here a detailed proof.

\begin{lemme}
\label{exo : 4.1}
Let $\Omega \subset \R^{N}$ be an open set. If $H_{t} = \{x \in \R^{N} | x \cdot e < t\} $ for
some $e \in \mathbb{S}^{N-1}$, $t \in \R$, and $E$ is a set of locally finite perimeter such that $|E| < +\infty$, then $E \cap H_{t}$ is a
set of finite perimeter in $\Omega$, and, for a.e. $t \in \R$,
\begin{equation}
\label{exo 13 :1}
    \mu_{E \cap H_{t}} = \mu_{E} \mres H_{t} + e \mathcal{H}^{N-1}\mres (E \cap \partial H_{t}) .
\end{equation}
Moreover, for a.e. $t \in \R$, we have
\begin{equation}
\label{exo : 13.2}
P(E \cap H_{t},\Omega ) \leqslant P(E,\Omega).
\end{equation}
If we further assume that $| E \cap \{x\cdot e > t\}\cap \Omega| > 0$, then inequality \eqref{exo : 13.2} is strict. 

\end{lemme}

\begin{rmq}
In particular, for $\Omega = \R^{N}$, Lemma~\ref{exo : 4.1} implies 
\begin{equation}
P(E \cap H_{t}) \leqslant P(E),
\end{equation}
with strict inequality if $| E \cap \{x\cdot e > t\}| > 0$.
\end{rmq}

\begin{proof}
The identity \eqref{exo 13 :1} is obtained by reasoning exactly as in \cite[Lemma 15.12]{maggi_sets_2012}. Thus, we only focus on the proof of inequality \eqref{exo : 13.2}. 

First, taking the total variation of \eqref{exo 13 :1} in $\Omega$ yields 
\begin{equation}
\label{decompo perimetre}
    P(E \cap H_{t},\Omega) = P(E,H_{t}\cap \Omega) + \mathcal{H}^{N-1}(E\cap \partial H_{t}\cap \Omega),
\end{equation}
because the supports of the measures $\mu_E \mres H_t$ and $\Haus \mres (E\cap \partial H_t)$ are disjoint.
Then, denoting $H_{t}^{+} := \{x \in \R^{N} | x \cdot e > t\}$, we notice that
\[P(E,\Omega) \geqslant P(E,H_{t}\cap \Omega) + P(E,H_{t}^{+}\cap \Omega).\]
Additionaly, for a set of finite perimeter $F$ such that $ |F| < \infty$, the proof of \cite[Proposition 19.22]{maggi_sets_2012} yields
 \[P(F,\partial H) \leqslant P(F,H) - \int_{H\cap \partial^* F} (1-|\nu_F \cdot e_N|) d \Haus \leqslant P(F,H).\]
 Applying this inequality with $F = E \cap H$, where $H$ is the superior half space $H=\{x_N > 0\}$, leads to
 \[\mathcal{H}^{N-1}(E^{(1)}\cap \partial H) \leqslant P(E\cap H, \partial H) \leqslant P(E \cap H,H) = P(E,H).\]
  Besides, \cite[Proposition 19.22]{maggi_sets_2012} yields strict inequality if we further assume that $|E\cap H|>0$.
 The proof of this proposition can be easily generalized to show that, for  every set of finite perimeter $E$ in $\Omega$ such that $|E| < \infty$, we have
 
 \[\mathcal{H}^{N-1}(E^{(1)}\cap \partial H_{t}^{+}\cap \Omega) \leqslant P(E,H_{t}^{+}\cap \Omega),\]
 with strict inequality if we further assume that $|E\cap H_{t}^{+} \cap \Omega|>0$.

Now, we claim that 
\[\mathcal{H}^{N-1}(E^{(1)}\cap \partial H_{t}^{+}\cap  \Omega) = \mathcal{H}^{N-1}(E\cap \partial H_{t}^{+}\cap  \Omega).\]
We recall that, for any measurable set $E$ such that $|E| < + \infty$, Lebesgue's density Theorem yields \[|E\setminus E^{(1)}|=0.\]

Hence, $|(E\cap \Omega )\setminus (E^{(1)}\cap \Omega)| = |E\setminus E^{(1)} \cap \Omega| \leqslant |E\setminus E^{(1)}| = 
0.$
Then, the coarea formula leads to
\[0 = |(E\cap \Omega )\setminus (E^{(1)}\cap \Omega)| = \int_{\R}\Haus(((E\cap \Omega )\setminus (E^{(1)}\cap \Omega))\cap \{x\cdot e =t\})dt.\]
Thus, for a.e. $t \in \R$, 
\[\mathcal{H}^{N-1}((E\cap \Omega )\setminus (E^{(1)}\cap \Omega)\cap \partial H_{t}^{+})  = 0.\]
This concludes the proof of our claim. Therefore, 
\begin{equation}
\label{inclusion frontiere}
\mathcal{H}^{N-1}(E\cap \partial H_{t}^{+}\cap \Omega) = \mathcal{H}^{N-1}(E^{(1)}\cap \partial H_{t}^{+}\cap \Omega) \leqslant P(E,H_{t}^{+}\cap \Omega)).
\end{equation}

Now, pluging \eqref{inclusion frontiere} in \eqref{decompo perimetre} yields
\begin{align*}
    P(E,\Omega) &\geqslant P(E,H_{t}\cap \Omega) + P(E,H_{t}^{+}\cap \Omega) \\
    &\geqslant P(E,H_{t}\cap \Omega) +\mathcal{H}^{N-1}(E\cap \partial H_{t}^{+}\cap \Omega) \\
    &= P(E \cap H_{t}, \Omega).
\end{align*}
Once again, the inequality is strict if we assume $|E\cap H_{t}^{+}\cap \Omega|>0$.
\end{proof}

\begin{rmq}
    Notice that this implies the same result for the half-space defined as $H_{t}^+ = \{x \in \R^{N} | x \cdot e > t\} $ because of the identity 
    \[H_{t}^+ = \{x \in \R^{N} | x \cdot (-e) < -t\}. \]
\end{rmq}

\begin{lemme}
\label{Maggi_cvx}
  Let $\Omega \subset \R^{N}$ be an open set and $E$ be a set of finite perimeter in $\Omega$ such that $|E| < + \infty $. If $K$ is a convex set, then $E \cap K$ is a set of finite
perimeter  in $\Omega$ and 
\[P(E \cap K, \Omega) \leqslant P(E,\Omega).\]
\end{lemme}

\begin{proof}
    Let $K$ be a convex set, therefore (see \cite[paragraph 7.5.3]{borwein2010convex} for instance) there exists countably many closed half spaces $\overline{H_{t_{n},e_{n}}} := \{ x \in \R^{d} | x\cdot e_{n} \leqslant t_{n} \},$ with $t_{n} \in \R$ and $e_{n} \in \mathbb S^{d-1}$, such that 
    \[\overline{K} = \bigcap_{n \in \N} \overline{H_{t_{n},e_{n}}}.\]
    The main idea of this proof is to explain how this Lemma follows from Lemma~\ref{exo : 4.1}. From \eqref{exo : 13.2}, except for a negligible set, denoted $\mathcal{N}$, intersecting with an open half set reduces the perimeter. 
    \begin{itemize}
        \item Let $t \in \R$ and $e\in \mathbb{S}^{d-1}$. Firstly, we recall that if $t \notin \mathcal{N}$ then \eqref{exo : 13.2} yields $P(E \cap \overline{H_{t,e}},\Omega) = P(E \cap{H_{t,e}},\Omega) \leqslant P(E,\Omega)$. Now, if $t \in \mathcal{N}$, since $\mathcal{N}$ is negligible, there exists a decreasing sequence $(\varepsilon_{n})$ converging to $0$ such that for all $n$, $t + \varepsilon_{n} \notin \mathcal{N}$ and thus
        \[P(E,\Omega) \geqslant P(E\cap H_{t+\varepsilon_{n},e},\Omega).\]

Letting $n$ tend to $0$, since $(H_{t+\varepsilon_{n},e})\cap E$ is a decreasing sequence, which converges to $H_{t,e}\cap E$ in $L^{1}$, by lower-semi-continuity of the perimeter this leads to

       \[P(E,\Omega) \geqslant P(E\cap H_{t,e},\Omega) = P(E\cap \overline{H_{t,e}},\Omega).\]
       \item Then, denoting $K_{n} = \bigcap_{k=0}^{n} \overline{H_{t_{k},e_{k}}}$, through recursion it follows from the argument above that $E\cap K_{n}$ is of finite perimeter in $\Omega$ and 
       \[P(E\cap K_{n},\Omega) \leqslant P(E,\Omega).\]
       Besides, from the monotone convergence theorem, $(K_{n}\cap E)$ converges towards $\overline{K}\cap E$ in $L^{1}(\R^{d})$. Therefore, the lower semi-continuity of the perimeter implies that $E\cap K$ is of finite perimeter and
       \[P(E\cap K,\Omega) = P(E\cap \overline{K},\Omega) \leqslant P(E,\Omega).\]
    \end{itemize} 
    \end{proof}
    
\subsection{Existence of solutions to Plateau's problem}
\label{Section existence}

In this short Section, we prove the existence of solutions for the problem \eqref{def : pblimit}. In particular this justifies the existence of Plateau-quasi-minimizers. 

\begin{prop}
    \label{existence}
    Plateau's problem \eqref{def : pblimit} admits a solution.
\end{prop}

\begin{proof}
  
    Firstly, since $E_0$ is a competitor, in the sense of Definition~\ref{competitors}, of finite perimeter in $\mathcal{D}$, it yields that the infimum 
	\[I := \inf \{ P(\Omega ,\mathcal{D}) ,\ \Omega \text{ is a competitor as defined in Definition~\ref{competitors}}\}\] is finite.

	Now, let $(\Omega_{n})_{n\in \N}$ be a minimizing sequence for problem~\eqref{def : pblimit}. Since $I$ is finite, one can assume that there exists a constant $C>0$ such that $\sup_{n\in \N} P(\Omega_{n},\mathcal{C}') \leqslant C$.
	By Rellich theorem, there exist a set $\Omega_{0}$ of finite perimeter in $\mathcal{D}$ and a subsequence $(\Omega_{n_{k}})$, such that $\chi_{\Omega_{n_{k}}} \to \chi_{\Omega_{0}}$ in $L^{1}(\mathcal{D})$. Furthermore, $\Omega_n$ are competitors, which leads them to satisfy the boundary condition $\Omega_{n} \setminus \overline{\mathcal{C}} = E_0 \setminus \overline{\mathcal{C}}$. Hence, this condition is also satisfied by $\Omega_0$. Since $\Omega_0$ is defined up to a negligible set, we can assume that $\Omega_0$ is a Borel set, making $\Omega_0$ a competitor.
	Finally, we know that the perimeter, in the De Giorgi sense, is lower semi-continuous. Thus, we get the desired inequality $P(\Omega_{0},\mathcal{C}') \leqslant \liminf P(\Omega_{n_{k}}, \mathcal{C}') = I$.
\end{proof}

\section{Ahlfors regularity up to the boundary}
\label{section AR}

\subsection{Ahlfors regularity result}

We start by proving the Ahlfors regularity of quasi-optimal surfaces up to the boundary. 

\begin{thm}[Ahlfors regularity]
\label{th : AhlReg}
     Let $(\mathcal{D}, \mathcal{C}, E_0)$ satisfying Hypothesis H (see Definition~\ref{def H}) and let $\Omega_{0}$ be a Plateau-quasi-minimizer satisfying $\spt \mu_{\Omega_0}=\partial \Omega_0$. Then, there exist some constants $C_{1},C_{2} > 0$, depending only on the dimension $N$ and the Ahlfors regularity constant of $\Sigma$, such that for all $x \in \partial \Omega_{0}\cap \mathcal{D} $ and $r>0$ such that $B(x,r)\subset \mathcal{D} $,
    \begin{equation}
    \label{eq : ARess2}
        C_{1} \leqslant \frac{P(\Omega_{0},B(x,r))}{r^{N-1}} \leqslant C_{2}.
    \end{equation}
    
In particular, $\Haus$-a.e. $\partial \Omega_0 \cap \mathcal{D} = \partial^* \Omega_0 \cap \mathcal{D}$.
   
\end{thm}

Notice that the only assumption on the boundary constraint $E_0$ needed for this result is hypothesis~\ref{H2}.

\subsection{Interior Ahlfors regularity}

Inside $\mathcal{C}$, the proof of the Ahlfors regularity is standard, see for instance \cite[Theorem 16.14]{maggi_sets_2012}. However, in order to extend the proof up to the boundary, we recall first the interior argument. 

\begin{prop}[Interior Ahlfors regularity]
\label{th : AhlRegInt}
    Let $(\mathcal{D}, \mathcal{C}, E_0)$ satsifying Hypothesis H (see Definition~\ref{def H}) and let $\Omega_{0}$ be a Plateau-quasi-minimizer, with constant $Q \geqslant 1$, satisfying $\spt \mu_{\Omega_0}=\partial \Omega_0$. Then there exist some constants, $C_{1},C_{2} > 0$ such that for all $x \in \partial \Omega_{0}\cap \mathcal{C}$ and $r>0$ satisfying $B(x,r) \subset \mathcal{C}$,
    \begin{equation}
    \label{ARint}
        C_{1} \leqslant \frac{P(\Omega_{0},B(x,r))}{r^{N-1}} \leqslant C_{2}.
    \end{equation}
    
    In particular, $\partial \Omega_0 \cap \mathcal{C} = \partial^* \Omega_0 \cap \mathcal{C}$.
\end{prop}

\begin{figure}[h]
\centering
 \includegraphics[scale = 0.25]{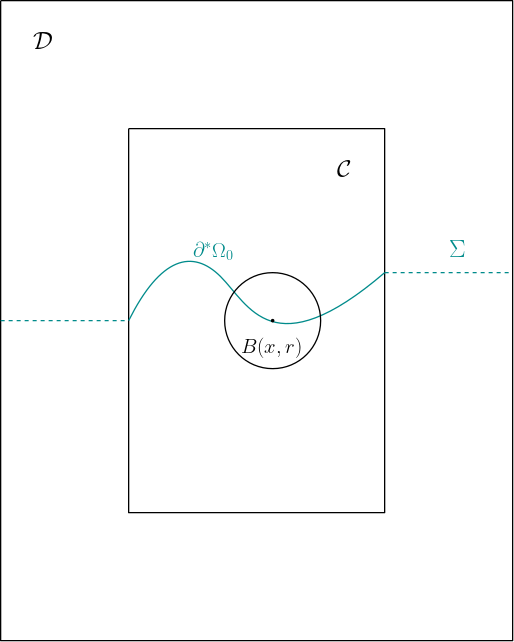}
  \label{fig : AR int}
    \caption{Interior Ahlfors regularity: $B(x,r) \subset \mathcal{C}$}
\end{figure}

\begin{proof}Let $x \in \partial \Omega_{0}\cap \mathcal{C}$. Since we made the assumption $\spt \mu_{\Omega_0}=\partial \Omega_0$,
\cite[Proposition 12.19]{maggi_sets_2012} yields 
    \begin{equation}
    \label{step 0}
        0 < |\Omega_{0}\cap B(x,r)|< \omega_{N}r^{N},
    \end{equation}
and thanks to \cite[Proposition 2.16]{maggi_sets_2012}, for almost every $r >0$, 
    \begin{equation}
    \label{th 2.2 : a.e}
         \Haus(\partial^{*}\Omega_{0} \cap \partial B(x,r))  = 0.
    \end{equation}
    Let us fix such $r$.
    
\textit{Step 1.} First, we select a competitor for Plateau's problem, so that we can use the quasi-minimality assumption. We consider 
\[\Omega = \Omega_{0} \setminus \overline{B(x,r)}.\]
 Since the changes are only made inside of $\overline{\mathcal{C}}$ (from the assumption $\overline{B(x,r)} \subset \overline{\mathcal{C}}$), $\Omega$ is a competitor for \eqref{def : pblimit}. And this yields 
    \[\Haus(\partial^* \Omega_0 \setminus \partial^* \Omega \cap \mathcal{D})\leqslant Q \Haus(\partial^* \Omega \setminus \partial^* \Omega_0 \cap \mathcal{D}).\]
    Let $s >r $ be such that $B(x,s) \subset \mathcal{D}$. Once again, given that modifications are only made in $\overline{B(x,r)} \subset B(x,s)$ and $Q\geqslant 1$, Remark~\ref{rq : quasimin} leads to
    \[P(\Omega_{0},B(x,s)) \leqslant Q P(\Omega,B(x,s)).\]
     We notice that $\overline{B(x,r)}$ is a smooth set and therefore the reduced boundary, essential boundary and topological boundary coincide. Moreover, the set of points of density $0$ of the ball is $(\overline{B(x,r)})^{c}$. Thus, \eqref{Maggi prive} yields
    \begin{align*}
    P(\Omega, B(x,s)) = &P(\Omega_{0},B(x,s) \setminus \overline{B(x,r)}) + P(B(x,r),\Omega_{0}^{(1)} \cap B(x,s))\\ &+ \Haus(\{\nu_{\Omega_{0}} = - \nu_{B(x,r)}\}\cap B(x,s)).
    \end{align*}
    But, by definition, 
    \begin{align*}
    \{\nu_{\Omega_{0}} = - \nu_{B(x,r)}\} &:= \{ x \in \partial_{*}\Omega_{0} \cap \partial B(x,r) \ | \ \nu_{\Omega_{0}}(x) = - \nu_{B(x,r)}(x)\} \\
    &\subset  \partial_{*}\Omega_{0} \cap \partial B(x,r),
\end{align*}     
    and we have chosen $r$ so that (by Federer's Theorem) the set $\partial_{*}\Omega_{0} \cap \partial B(x,r)$ is $\Haus$-negligible. Besides, since $s>r$,
    \begin{align*}
    P(B(x,r),\Omega_{0}^{(1)} \cap B(x,s)) &= \Haus(\partial B(x,r) \cap \Omega_{0}^{(1)}\cap B(x,s)) \\
    &= \Haus(\partial B(x,r) \cap \Omega_{0}^{(1)}).
    \end{align*} 
     Hence,
    \begin{align*}
    P(\Omega_{0},B(x,s)) &\leqslant Q P(\Omega,B(x,s)) \\
    &= Q P(\Omega_{0},B(x,s) \setminus \overline{B(x,r)}) + Q \Haus(\partial B(x,r) \cap \Omega_{0}^{(1)}).
\end{align*}     
    And then by letting $s$ go to $r$ we get
    \begin{equation}
    \label{step 1}
        P(\Omega_{0},B(x,r)) \leqslant Q \Haus(\partial B(x,r) \cap \Omega_{0}^{(1)}),
    \end{equation}
    which clearly implies the upper bound of \eqref{ARint}, with constant only depending on $Q$ and $N$.
    
    \textit{Step 2.} The observation made above with \eqref{th 2.2 : a.e} and \eqref{Maggi inter} yield 
    \[P(\Omega_{0} \cap B(x,r)) = P(\Omega_{0},B(x,r)) +\Haus(\Omega_{0}^{(1)} \cap \partial B(x,r)).\]
    So by adding $\Haus(\Omega_{0}^{(1)} \cap \partial B(x,r))$ to both sides of \eqref{step 1}, we get 
    \[P(\Omega_{0} \cap B(x,r)) \leqslant (1 +Q) \Haus(\Omega_{0}^{(1)} \cap \partial B(x,r)).\]
    Next, applying the isoperimetric inequality (see \cite[Theorem 14.1]{maggi_sets_2012}) to $\Omega_{0} \cap B(x,r)$ 
    yields
    \begin{equation}
        \label{step 2}
        N \omega_{N}^{1/N}|\Omega_{0} \cap B(x,r)|^{\frac{N-1}{N}} \leqslant P(\Omega_{0} \cap B(x,r)) \leqslant (1+Q)\Haus(\Omega_{0}^{(1)} \cap \partial B(x,r)).
    \end{equation}\\   
    \textit{Step 3.} We consider the function $m : s \mapsto |\Omega_{0}\cap B(x,s)|= |\Omega_{0}^{(1)}\cap B(x,s)| = \int_0^s \Haus(\Omega_0 \cap B(x,r))ds$. From the coarea formula, we know that $m$ is differentiable almost everywhere and $m'(s) = \Haus(\Omega_{0}^{(1)}\cap \partial B(x,s))$ for a.e. $s>0$. Then by \eqref{step 2} we have :
    \begin{equation}
        \label{step 3}
        N\omega_{N}^{1/N}m(s)^{\frac{N-1}{N}} \leqslant (1+Q)m'(s).
    \end{equation}
    But, thanks to \eqref{step 0}, we know that $m$ is positive and also that $m(0^{+})=0$. So by dividing by $m(s)^{\frac{N-1}{N}}$, integrating \eqref{step 3} on $]0,r[$ and setting it to the power of $N$, we get
    \begin{equation*}
        w_{N}r^{N} \leqslant (1+Q)^{N}m(r).
    \end{equation*}
    Applying the same reasoning to $\R^{N}\setminus \Omega_{0}$ instead of $ \Omega_{0}$ we get 
    \[|B(x,r)\setminus  \Omega_{0}| \geqslant \frac{\omega_{N}r^{N}}{(1+Q)^{N}}.\]
    Hence, from the decomposition $\omega_N r^N = |B(x,r)| = |B(x,r) \cap \Omega_0| + |B(x,r)\setminus \Omega_0|$,
    \begin{equation}
        \label{step 3.0}
        \frac{1}{(1+Q)^{N}} \leqslant \frac{|\Omega_{0}\cap B(x,r)|}{\omega_{N}r^{N}}\leqslant 1-\frac{1}{(1+Q)^{N}}.
    \end{equation}
    
The inequality~\ref{step 3.0} yields in particular that $x$ is in the essential boundary of $\Omega_0$ and therefore, $\partial \Omega_0 \cap \mathcal{C} = \partial^* \Omega_0 \cap \mathcal{C} $. 
    
    \textit{Step 4.} We now apply the relative isoperimetric inequality, see for instance \cite[Proposition 12.37]{maggi_sets_2012} with $t=1-\frac{1}{(1+Q)^{N}}$, and we obtain the existence of a constant $\gamma_{N}>0$ such that 
    \begin{align*}
    P(\Omega_{0},B(x,r)) &\geqslant \gamma_{N}|\Omega_{0}\cap B(x,r)|^{\frac{N-1}{N}} \geqslant \gamma_{N}(\frac{1}{(1+Q)^{N}}\omega_{N}r^{N})^{\frac{N-1}{N}}\\ 
    &=\gamma_{N}\frac{1}{(1+Q)^{N-1}}\omega_{N}^{\frac{N-1}{N}}r^{N-1}.
    \end{align*}
    This concludes the proof of the interior regularity for almost every $r$, with \[C_{1} = \gamma_{N}\frac{1}{(1+Q)^{N-1}}\omega_{N}^{\frac{N-1}{N}}.\]    
    
    \textit{Step 5.} Let us explain how to extend the previous result to any $r>0$. Take $r>0$ belonging to the negligible set that does not satisfy \eqref{th 2.2 : a.e}. Since this set in negligible, there exist $\frac{r}{2}<t<r$ satisfying \eqref{th 2.2 : a.e}. Thus,
    \[\frac{C_{1}}{2^{N-1}}r^{N-1} \leqslant C_{1}t^{N-1} \leqslant P(\Omega_{0},B(x,t)) \leqslant P(\Omega_{0},B(x,r)).\]
    Similarly, there exist $r < t < 2r$ such that $B(x,t) \subset \mathcal{C}$ satisfying \eqref{th 2.2 : a.e}. Thus,
    \[2^{N-1}C_{2} r^{N-1} \geqslant C_{2}t^{N-1} \geqslant P(\Omega_{0},B(x,t)) \geqslant P(\Omega_{0},B(x,r)).\]
    This concludes the proof of the interior regularity for all $r$.  
\end{proof}

\subsection{Ahlfors regularity at the boundary}

Now, we establish the Ahlfors regularity at the boundary of Plateau-quasi-minimizers, by using either the interior Ahlfors regularity or the assumption of Ahlfors regularity of $\Sigma$.

\begin{prop}[Ahlfors regularity at the boundary]
\label{th : AhlReg bdd}
    Let $(\mathcal{D}, \mathcal{C}, E_0)$ satisfying Hypothesis H (see Definition~\ref{def H}) and let $\Omega_{0}$ be a Plateau-quasi-minimizer (see Definition~\ref{def quasi min}) with constant $Q>0$, satisfying $\spt \mu_{\Omega_0}=\partial \Omega_0$. Then, there exist some constants, $C_{1},C_{2} > 0$ such that for all $x \in \partial \Omega_{0}\cap \mathcal{D}$ and $r>0$ such that $B(x,r) \subset \mathcal{D}$ and $B(x,r) \cap \partial \mathcal{C} \neq \emptyset$,
    \begin{equation}
    \label{eq : AR bdd}
        C_{1} \leqslant \frac{P(\Omega_{0},B(x,r))}{r^{N-1}} \leqslant C_{2}.
    \end{equation}
\end{prop}

\begin{proof}

\begin{figure}[ht]
\centering
 \includegraphics[scale = 0.25]{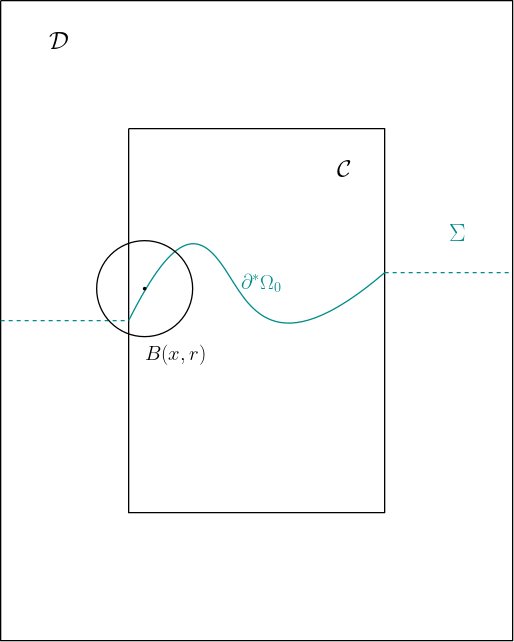}
  \label{fig : AR bound}
    \caption{Ahlfors regularity at the boundary: $B(x,r) \cap \partial \mathcal{C} \neq \emptyset$}
\end{figure}
 
 First we explain the argument to get the upper bound of \eqref{eq : AR bdd}, then we will focus on the lower bound, for each we distinguish several possible scenarios.

     \textit{The upper bound.} Let $x \in \partial \Omega_{0}\cap \mathcal{D}$ and $r>0$ such that $B(x,r) \subset \mathcal{D}$. We distinguish three possible cases. 
     \begin{itemize}
     \item If $x \in \mathcal{D}\setminus \overline{\mathcal{C}}$ then the boundary constraint yields that $x\in \Sigma $. Thus, we consider $\Omega' := ((E_0 \setminus \overline{\mathcal{C}})\cap B(x,r)) \cup (\Omega_0 \setminus B(x,r))$, which is a competitor, and the quasi-minimality of $\Omega_0$  and Remark~\ref{rq : quasimin} imply
    \[P(\Omega_0,B(x,r)) \leqslant Q P(\Omega',B(x,r)).\]
     Furthermore, recall that the topological boundary of a convex set is Ahlfors regular and by assumption~\ref{H2} $\Sigma$ is also Ahlfors regular. Hence, 
\begin{align*}
      P(\Omega',B(x,r)) &\leqslant \Haus(\Sigma \cap B(x,r)) + \Haus(\partial \mathcal{C}\cap B(x,r)) \leqslant  C_2 r^{N-1},
\end{align*}     
     and the upper bound is proved in this case.
     
     \item Otherwise, if $x \in \overline{\mathcal{C}}$ and $B(x,r) \cap \Sigma \neq \emptyset$, then we set $y \in B(x,r) \cap \Sigma$. Hence, $B(x,r) \subset B(y,2r)$ and as in the first case, we consider the competitor $\Omega' := ((E_0 \setminus \overline{\mathcal{C}})\cap B(y,2r)) \cup (\Omega_0 \setminus B(y,2r))$. Then, the Ahlfors regularity of $\Sigma$ and $\partial \mathcal{C}$ implies
     \[P(\Omega_0,B(x,r)) \leqslant P(\Omega_0,B(y,2r)) \leqslant Q P(\Omega',B(y,2r)) \leqslant Q C_2 2^{N-1} r^{N-1}.\]
     \item Finally, assume that $x \in \overline{\mathcal{C}}$ and $B(x,r) \cap \Sigma = \emptyset$. In this case, we use the same strategy as in the proof of the interior Ahlfors regularity. Compared to the interior case, here $\Omega = \Omega_0 \setminus \overline{B(x,r)}$ is no longer necessarily a competitor for Plateau's problem. More precisely, in order to be able to use the same argument as in the interior case, we need to ensure that for $s>r$,
     \begin{equation}
    \label{lem 3.4}
    P(\Omega_0,B(x,s)) \leqslant Q P(\Omega,B(x,s)).
    \end{equation}
Once this inequality is proved the rest of the proof done for the interior regularity unfolds the same way to get the estimate. In particular it will also imply that $x \in \partial^*\Omega_0$. There remains thus only to show \eqref{lem 3.4}. This is done in Lemma~\ref{lem AR}.

     \end{itemize}

\textit{The lower bound.} We will now focus on the lower bound, which, as we will show, derives from the assumption that $\Sigma$ is Ahlfors regular and the interior Ahlfors regularity. As mentioned earlier, we distinguish two cases: either a substantial part of the surface coincides with $\Sigma$ and we use the Ahlfors regularity of $\Sigma$, or we apply the interior regularity. Let $x \in \partial \Omega_{0}\cap \mathcal{D}.$ 
\begin{enumerate}
    \item First, if $d(x,\Sigma) \leqslant \frac{r}{2}$ there exists a point $y \in \Sigma $ such that $|x-y| \leqslant \frac{r}{2}$. Thus, since by boundary constraint $\partial^*\Omega_0 \cap B(y,\frac{r}{2}) \supset \partial^*E_0 \cap B(y,\frac{r}{2})\setminus \overline{\mathcal{C}} = \Sigma \cap B(y,\frac{r}{2})$,
   the Ahlfors regularity of $\Sigma$ yields 
    \[\frac{C_1}{2^{N-1}} r^{N-1} \leqslant \Haus(\Sigma \cap B(y,\frac{r}{2}))\leqslant P(\Omega_0,B(y,\frac{r}{2})) \leqslant  P(\Omega_0,B(x,r)).\]
        
    \item Else, if $d(x,\Sigma) > \frac{r}{2}$, the interior regularity (Proposition~\ref{ARint}) will allow us to conclude. In this case, $B(x, \frac{r}{2})\cap \Sigma = \emptyset$ thus the same strategy as in the last case for the upper bound is developped. Namely, we use Lemma~\ref{lem AR} to justify that the same proof as in the interior case can be applied. 
    \end{enumerate}  

\end{proof}

To complete the proof at the boundary, we now need to prove the key argument for the lower bound, in the second case, i.e. when $d(x, \Sigma) > \frac{r}{2}$. This is done in the following Lemma.

\begin{lemme}
\label{lem AR}
Let $(\mathcal{D}, \mathcal{C}, E_0)$ satisfying Hypothesis H (see Definition~\ref{def H}) and let $\Omega_{0}$ be a Plateau-quasi-minimizer. Then, for every $x \in \partial \Omega_{0}\cap \mathcal{D} $ and almost every $r>0$ such that $ B(x,r) \subset \mathcal{D}$ and $B(x,r)\cap \Sigma = \emptyset$, we have for all $s>r$
\begin{equation}  
        P(\Omega_{0},B(x,s))\leqslant Q P(\Omega_{0}\setminus \overline{B(x,r)},B(x,s)).
    \end{equation}

    \end{lemme}
    
\begin{proof}
To prove this Lemma, we will select a proper competitor and use the quasi-minimality of $\Omega_0$. 

First, \cite[Proposition 2.16]{maggi_sets_2012} claims that for almost every $r >0$, 
    \begin{equation}
    \label{lem 3.5 : a.e}
         \Haus(\partial^{*}\Omega_{0} \cap \partial B(x,r))= 0.
    \end{equation}
    Let us fix such $r$. We distinguish two possible cases, as summarized in Figure~\ref{fig : lem3.5-2}.

\begin{figure}[ht]
\centering
\subfloat[\centering Case 1.]{\includegraphics[width = 0.4 \textwidth]{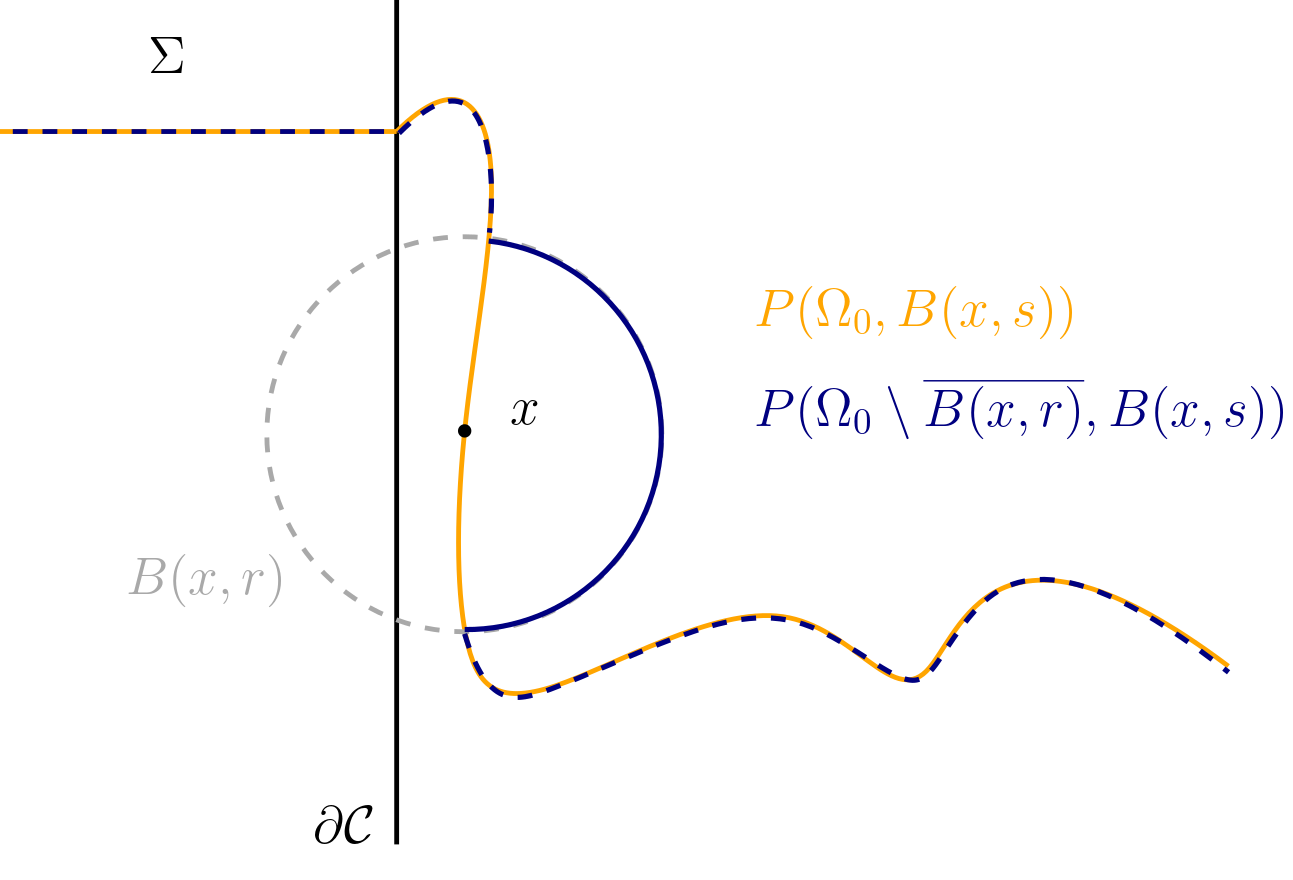} \label{fig : AR bord 1}}
\quad
\subfloat[\centering Case 2.]{\includegraphics[width = 0.50 \textwidth]{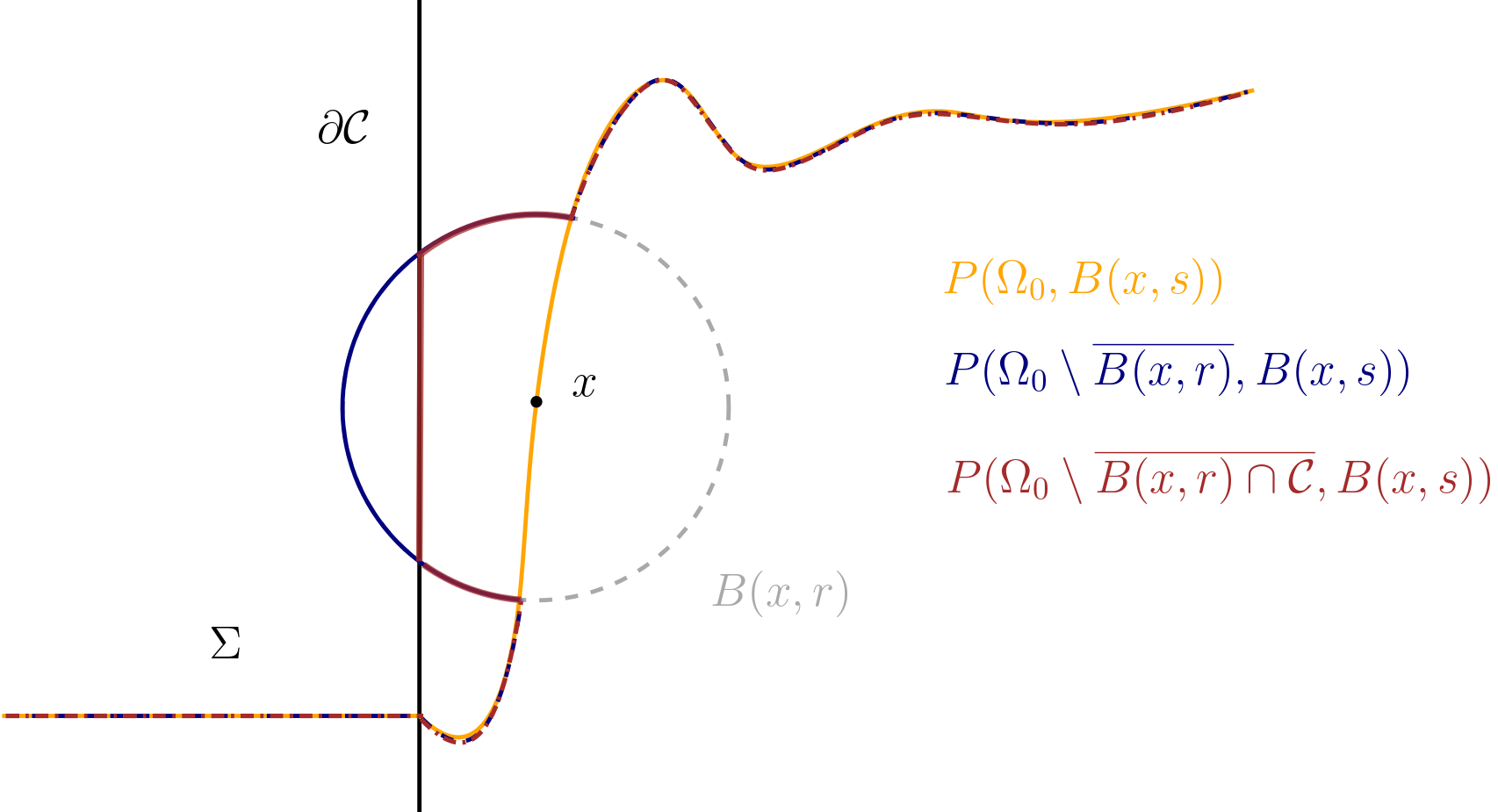} \label{fig : AR bord 2}}
    \caption{Lemma 3.5}
    \label{fig : lem3.5-2}
\end{figure}

\textit{Case 1.} If $\Omega_0 \cap \overline{B(x,r)} \subset \overline{\mathcal{C}}$ (Figure~\ref{fig : AR bord 1}), in which case, since $\Omega_0 \setminus \overline{B(x,r)}$ is a competitor, the quasi-minimality of $\Omega_0$ and Remark~\ref{rq : quasimin} are enough to conclude. 

\textit{Case 2.} Otherwise if $\Omega_0 \cap \overline{B(x,r)} $ is not contained in $\overline{\mathcal{C}}$, then $\Omega_0 \setminus \overline{B(x,r)}$ is no longer a competitor (Figure~\ref{fig : AR bord 2}). However, $\Omega_0 \setminus \overline{(B(x,r)\cap \mathcal{C}})$ is a competitor for \eqref{def : pblimit}, hence Remark~\ref{rq : quasimin} yields
\begin{equation}
P(\Omega_0, B(x,s)) \leqslant QP(\Omega_0 \setminus \overline{B(x,r)\cap \mathcal{C}},B(x,s)).
\end{equation}
Besides, thanks to \eqref{Maggi prive}, we can rewrite 
\begin{align*}
P(\Omega_0 \setminus \overline{B(x,r)\cap \mathcal{C}},B(x,s)) &= P(\Omega_0,B(x,s) \setminus \overline{B(x,r)\cap \mathcal{C}})\\
 & \quad + P(\overline{B(x,r) \cap \mathcal{C}}, \Omega_0^{(1)}\cap B(x,s))\\
 & \quad + \Haus(\{\nu_{\Omega_0} = - \nu_{B(x,r)\cap \mathcal{C}}\})\\
&=P(\Omega_0,B(x,s) \setminus \overline{B(x,r)\cap \mathcal{C}}) \\
& \quad + P(B(x,r) \cap \mathcal{C}, \Omega_0^{(1)}\cap B(x,s))\\
 & \quad + \Haus(\{\nu_{\Omega_0} = - \nu_{B(x,r)\cap \mathcal{C}}\}),
\end{align*}
and 
\begin{align*}
P(\Omega_0 \setminus \overline{B(x,r)},B(x,s)) &= P(\Omega_0,B(x,s) \setminus \overline{B(x,r)})
  + P(\overline{B(x,r) }, \Omega_0^{(1)}\cap B(x,s))\\
&=P(\Omega_0,B(x,s) \setminus \overline{B(x,r)}) + P(B(x,r) , \Omega_0^{(1)}\cap B(x,s)),
\end{align*}
because of \eqref{lem 3.5 : a.e}. Moreover, the assumption that $\Sigma\cap B(x,r) = \emptyset$, i.e. $\partial ^{*}\Omega_0 \cap B(x,r) \subset \overline{\mathcal{C}}$ ensures the equality 
\[\partial ^{*}\Omega_0 \cap (B(x,s) \setminus \overline{B(x,r)\cap \mathcal{C}}) =\partial ^{*}\Omega_0 \cap (B(x,s) \setminus \overline{B(x,r)}).\]
Thus 
\[P(\Omega_0,B(x,s)\setminus \overline{B(x,r)\cap \mathcal{C}}) = P(\Omega_0,B(x,s)\setminus \overline{B(x,r)}).\]
And then,
\begin{align*}
P(\Omega_0,B(x,s)) &\leqslant Q P(\Omega_0 \setminus \overline{B(x,r)\cap \mathcal{C}},B(x,s))\\
&= Q\left(P(\Omega_0,B(x,s) \setminus \overline{B(x,r)})+ P(B(x,r) \cap \mathcal{C}, \Omega_0^{(1)}\cap B(x,s))\right.\\
&\left. \quad  + \Haus(\{\nu_{\Omega_0} = - \nu_{B(x,r)\cap\mathcal{C}}\})\right)
\end{align*}
Notice that, by definition, $\Haus(\{\nu_{\Omega_0} = - \nu_{B(x,r)\cap\mathcal{C}}\}) \leqslant \Haus(\partial^{*}\Omega_0 \cap \partial(B(x,r)\cap \mathcal{C}))$. Next, we decompose : $B(x,s) \cap \Omega_0^{(1)} = B(x,s) \setminus (\partial^{*}\Omega_0 \cup \Omega_0^{(0)})$.

\begin{align*}
P(\Omega_0,B(x,s)) &\leqslant Q\left(P(\Omega_0,B(x,s) \setminus \overline{B(x,r)})+ \Haus(\partial^{*}\Omega_0 \cap \partial(B(x,r)\cap \mathcal{C}))\right.\\
& \quad +  P(B(x,r) \cap \mathcal{C}, B(x,s)) \\
&\left. \quad  - P(B(x,r)\cap \mathcal{C},B(x,s)\cap \Omega_0^{(0)}) -\Haus(\partial^* \Omega_0\cap \partial(B(x,r)\cap\mathcal{C}))\right).
\end{align*}
Finally, we use the convexity of $\mathcal{C}$ and Lemma~\ref{Maggi_cvx}. 
\begin{align*}
P(\Omega_0,B(x,s)) &\leqslant Q\left(P(\Omega_0,B(x,s) \setminus \overline{B(x,r)})+ P(B(x,r), B(x,s))\right.\\
&\left. \quad  - P(B(x,r)\cap \mathcal{C},B(x,s)\cap \Omega_0^{(0)})\right)\\
&= Q\left(P(\Omega_0,B(x,s) \setminus \overline{B(x,r)})+ P(B(x,r), B(x,s))\right.\\
&\left. \quad  - P(B(x,r),B(x,s)\cap \Omega_0^{(0)})\right)\\
&= Q\left(P(\Omega_0,B(x,s) \setminus \overline{B(x,r)})+ P(B(x,r), B(x,s)\cap \Omega_0^{(1)})\right)\\
&=QP(\Omega_0 \setminus \overline{B(x,r)},B(x,s)),
\end{align*}

thanks to \eqref{lem 3.5 : a.e}. The first inequality from the previous calculation is represented in  Figure \ref{fig : lem 3.5 2}.

\begin{figure}[ht]

\centering
\includegraphics[width = 0.5 \textwidth]{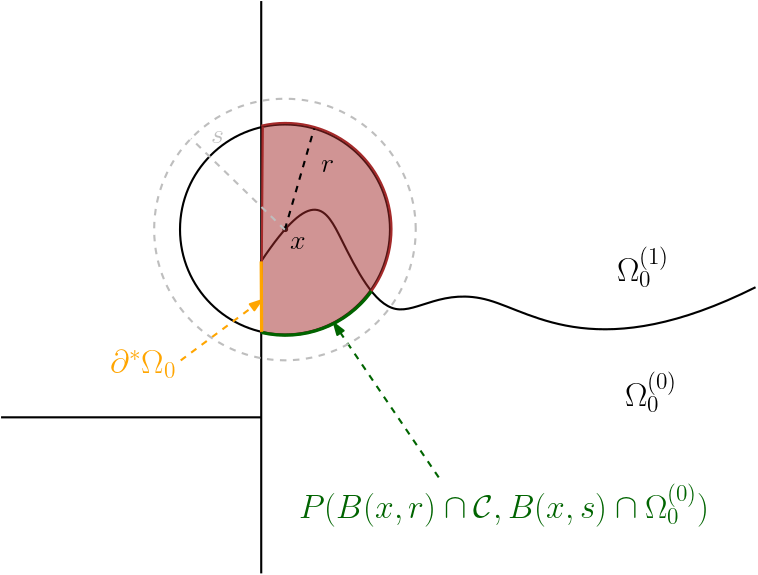} 

\caption{$P(B(x,r)\cap \mathcal{C},B(x,s)\cap \Omega_0^{(0)}) = P(B(x,r),B(x,s)\cap \Omega_0^{(0)})$} 
\label{fig : lem 3.5 2}
\end{figure}

This completes the proof in this case.
\end{proof}

This allow us to prove Theorem~\ref{th : AhlReg}.

\begin{proof}{\emph{(Theorem~\ref{th : AhlReg})}}
Let $x \in \partial \Omega_{0}\cap \mathcal{D}$ and $r>0$ be such that $B(x,r) \subset \mathcal{D} $. Then we are in one of the previously studied cases. If $B(x,r) \subset \mathcal{C}$ we apply the interior Ahlfors regularity, Proposition~\ref{ARint}. If $B(x,r) \subset \mathcal{D}\setminus \overline{ \mathcal{C}}$ we use the Alhfors regularity of $\Sigma$. Otherwise, we have $B(x,r) \cap \partial \mathcal{C} \neq \emptyset$ and we can use the Ahlfors regularity at the boundary, Proposition~\ref{th : AhlReg bdd}. We conclude by taking the infimum of the three lower bounds ($C_1$) and the maximum of the upper bounds ($C_2$).

Notice also that on the one hand the interior regularity leads to $\partial \Omega_0 \cap \mathcal{C} = \partial^* \Omega_0 \cap \mathcal{C}$. And on the other hand, the boundary constraint yields that $\partial^* \Omega_0 \cap \mathcal{D} \setminus \overline{\mathcal{C}} = \partial^* E_0 \cap \mathcal{D} \setminus \overline{\mathcal{C}} = \partial E_0 \cap \mathcal{D} \setminus \overline{\mathcal{C}} = \partial \Omega_0 \cap \mathcal{D} \setminus \overline{\mathcal{C}}$. Thus, there remains to justify why $\Haus$-a.e. $\partial \Omega_0 \cap \partial \mathcal{C} = \partial^* \Omega_0 \cap \partial \mathcal{C}$. In Proposition~\ref{th : AhlReg bdd}, we have shown that for $x \in \partial \Omega_0$ such that there exists a $r>0$ such that $B(x,r)\cap \Sigma = \emptyset$, $x \in \partial^*\Omega_0$. Therefore, since for $\Haus$-almost every $x\in \partial \Omega_0 \cap \partial \mathcal{C}$ there exists a $r>0$ such that $B(x,r) \cap \Sigma = \emptyset$, we can conclude that $\partial \Omega_0 \cap \mathcal{D} = \partial^* \Omega_0 \cap \mathcal{D}$, $\Haus$-a.e.
\end{proof}

%

\section{Uniform rectifiability up to the boundary}
\label{section unif rect}
\subsection{Definition and characterization of uniform rectifiability}

Having established the Ahlfors regularity of Plateau-quasi-minimizers up to the boundary, we will now prove that Plateau-quasi-minimizers satisfy Condition B, which yields in particular the uniform rectifiability of the quasi-optimal surface.

Condition B (see Definition~\ref{condition B}) is indeed a characterization of having Big Pieces of Lipschitz Graphs (or BPLG, for short), which implies in particular the  uniform rectifiability. Although the BPLG property is known to be strictly stronger than uniform rectifiability, the difference between the two is not significant (see \cite{david1993quantitative} for more information on the distinction between these two properties). For the curious reader, we first recall the definition of uniform rectifiability and having BPLG.

\begin{definition}[Uniform rectifiability]
Let $S \subset \mathcal{D}$ be an Ahlfors regular set in $\mathcal{D}$, see Definition~\ref{def AR}. We say that $S$ is uniformly rectifiable in $\mathcal{D}$, if there exists a constant $C>0$ such that for all $x\in S$ and all $r>0$, such that $B(x,r) \subset \mathcal{D}$, there exist a compact set $F \subset S\cap B(x,r)$ and a map $f : F \longrightarrow \R^{N-1}$ satisfying
\begin{equation}
\Haus (F) \geqslant C^{-1}r^{N-1},
\end{equation}
and for all $y,z \in F$
\begin{equation}
C^{-1}|z-y| \leqslant |f(z)-f(y)| \leqslant C|z-y|.
\end{equation}
\end{definition}

This definition means that a uniformly rectifiable set is locally the image of a bi-Lipschitz function up to a set of small $\Haus$-measure (with respect to $r$). This notion is investigated in detail in \cite{david1996uniform} and \cite{david2000uniform}. 

\begin{definition}[BPLG]
We say that $S$ has BPLG (Big Pieces of Lipschitz Graph), if $S$ is Ahlfors regular and if there exist constants $C>0$ and $\theta > 0$ such that for all $x \in S$ and $r>0$, there is an $(N-1)$-dimensional Lipschitz graph $G$ with constant smaller than $C$ such that 
\[\Haus(S\cap B(x,r) \cap G) \geqslant \theta r^{N-1}.\]
\end{definition}

\begin{prop}{\cite[Proposition 2, p. 103-104]{david1988morceaux}}
\label{prop : caract unif rect}
If a set $E$ satisfies Condition B in $\mathcal{D}$ (see Definition \ref{condition B}), then its boundary $\partial E$ has BPLG (Big Pieces Lipschitz Graph) in $\mathcal{D}$, in particular it is uniformly rectifiable in $\mathcal{D}$. 
\end{prop}

\subsection{Uniform rectifiability result}

\begin{thm}[Uniform rectifiability]
\label{th : unif rect}
    Let $(\mathcal{D}, \mathcal{C}, E_0)$ satisfying Hypothesis H (see Definition~\ref{def H}) and let $\Omega_{0}$ be a Plateau-quasi-minimizer such that $\spt \mu_{\Omega_0}=\partial \Omega_0$. Then, up to considering an equivalent set with the same assumption, $\partial \Omega_{0}$ has BPLG, so in particular is uniformly rectifiable. 
\end{thm}

Notice that only the assumptions~\ref{H2} and \ref{H3} on the boundary constraint $E_0$ are needed for this result.

To prove Theorem~\ref{th : unif rect}, thanks to Proposition~\ref{prop : caract unif rect}, we will prove that Condition B (see Definition~\ref{condition B}) is satisfied up to the boundary. Note that (thanks to Hypothesis~\ref{H2}) Theorem~\ref{th : AhlReg} insures the Ahlfors regularity, required in Condition B.

 As we did in Section~\ref{section AR}, the regularity outside and inside the convex set $\mathcal{C}$ will be carried out separately, and the results will then be extended up to the boundary.

\subsection{Interior uniform rectifiability}

First, we focus on justifying that Condition B is satisfied inside $\mathcal{C}$. 

\begin{prop}[Condition B inside $\mathcal{C}$]
\label{prop : rect int}
 Let $(\mathcal{D}, \mathcal{C}, E_0)$ satisfying  Hypothesis H (see Definition~\ref{def H}) and let $\Omega_{0}$ be a Plateau-quasi-minimizer such that $\spt \mu_{\Omega_0} = \partial \Omega_0$. Then, up to considering an equivalent set open and with the same assumption on the boundary, there exist a constant $C>1$ such that for all $x\in \partial \Omega_{0} \cap \mathcal{C}$ and $r>0$ such that $B(x,r) \subset \mathcal{C}$, there exist two balls $B_1 \subset B(x,r)\cap \Omega_0$ and $B_2 \subset B(x,r) \setminus \overline{\Omega_0}$ with radius larger than $C^{-1}r$. 
\end{prop}

This Proposition was first proved by David and Semmes in \cite{david1998quasiminimal}. However we will use a slightly different proof, detailed in \cite{rigot2000uniform}, as it is easily adaptable to adress the regularity at the boundary. The proof, based on several lemmas, consists in constructing a proper equivalent open set. But first, we recall a classical result useful for the following proof, the average formula for finite measure. 

\begin{lemme}
    \label{average formula}
    Let $\mu$ be a finite measure on $\R^N$ and $A$ a $\mu$-mesurable set. Let $f$ be a $\mu$ measurable function on $A$ such that $\int_{A} f d\mu \geqslant 0$. Then, there exists $t_{0} \in A$ such that 
    \begin{equation}
        \int_{A} f d\mu \geqslant \mu(A) f(t_{0}).
    \end{equation}
\end{lemme}

\begin{proof}
Assume, without loss of generality, that $\mu (A) > 0$.
We recall that, for all $t>0$,
\[\mu(\{f>t\}) = \int_{\{f>t\}}d\mu < \frac{1}{t}\int_{\{f>t\}}f d\mu \leqslant \frac{1}{t}\int_{A}f d\mu.\]

If $m = \frac{\int_{A}f d\mu}{\mu(A)} >0$ we apply this inequality with $t = m >0 $, and we get:
\[\mu(\{f>m\}) < \frac{1}{m}\int_{A}f d\mu = \mu(A).\]
Thus there must exist a $t_{0} \in A$ such that $f(t_{0}) \leqslant m = \frac{\int_{A}f d\mu}{\mu(A)}$. 
In other words
\[\int_{A} f d\mu \geqslant \mu(A) f(t_{0}).\]

 Otherwise, if $m=0$. Either $f$ is positive, then $f = 0$ $\mu$-a.e. Or, there exists a $t_{0}$ such that $f(t_{0}) \leqslant 0$. In any case, the desired inequality is satisfied.  
\end{proof}

\begin{lemme}
\label{lem 4.10}
There exists a constant $\varepsilon >0$ such that if $\Omega_0$ is a Plateau-quasi-minimizer then 
\[\spt \mu_{\Omega_0} \cap \mathcal{C}= \{ x\in \R^{N} : h(x,r) \geqslant \varepsilon, \forall r>0 \text{ such that } B(x,r) \subset \mathcal{C}\},\]
where $h(x,r) := r^{-N} \min \{|\Omega_0 \cap B(x,r)|,|B(x,r) \setminus \Omega_0|\}$.
\end{lemme}

\begin{proof}

 The proof of this Lemma is presented in \cite[Lemma 3.4]{rigot2000uniform}. For the reader's convenience, but mostly to underline the argument that needs to be adapted for the regularity at the boundary, we detail the proof of this crucial Lemma.

Let $\Omega_{0}$ be a Plateau-quasi-minimizer, and $x \in \spt \mu_{\Omega_0}\cap \mathcal{C}$ and $r>0$ be such that $B(x,r) \subset \mathcal{C}$. 

\textit{Step 1.} Assume first that $h(x,r) = r^{-N}|\Omega_0\cap B(x,r)|$. By the coarea formula and the average formula (see Lemma~\ref{average formula}), we can find $t \in ]r/2,r[$ such that 

\begin{align*}
r^{N}h(x,r) = |\Omega_0 \cap B(x,r)| = \mathcal{H}^N(\Omega_0 \cap B(x,r)) &= \int_0^r \Haus(\Omega_0 \cap \partial B(x,s))ds \\
&\geqslant \int_{r/2}^r \Haus(\Omega_0 \cap \partial B(x,s))ds \\
&\geqslant \frac{r}{2} \Haus(\Omega_0 \cap \partial B(x,t)).
\end{align*}
Thus, 
\begin{equation}
\label{lem 4.10 id1}
\Haus(\Omega_0 \cap \partial B(x,t)) \leqslant 2r^{N-1}h(x,r).
\end{equation}
Besides, we can choose $t$ such that $\Haus(\partial^{*}\Omega_0 \cap \partial B(x,t)) = 0$, since from \cite[Proposition 2.16]{maggi_sets_2012} almost every $t$ satisfies this property. Hence, in the average formula we can integrate over the set $]r/2,r[\setminus \mathcal{N}$, where $\mathcal{N}$ is the negligible set for which the previously mentioned equality does not hold. Thus, for $t'>t$, $\Haus(\{\nu_{\Omega_0} = -\nu_{\overline{B(x,t)}}\}\cap B(x,t'))=0$, and \eqref{Maggi prive} yields 
\begin{align}
\label{lem 4.10 id2}
P(\Omega_0 \setminus \overline{B(x,t)},B(x,t')) = P(\Omega_0, B(x,t')\setminus \overline{B(x,t)}) + \Haus(\Omega_0 \cap \partial B(x,t)).
\end{align}
Now, we use the fact that $\Omega_0$ is a P-quasi-minimizer and Remark~\ref{rq : quasimin}. We consider $\Omega_0 \setminus \overline{B(x,t)}$. Since $t<r$, this set is a competitor for \eqref{def : pblimit} therefore, combining \eqref{lem 4.10 id2} and \eqref{lem 4.10 id1}, we get the upper bound
\begin{align*}
P(\Omega_0, B(x,t')) &\leqslant Q P(\Omega_0 \setminus \overline{B(x,t)}, B(x,t')) \\
&= Q\left(P(\Omega_0, B(x,t')\setminus \overline{B(x,t)}) + \Haus(\Omega_0 \cap \partial B(x,t))\right)\\
&\leqslant Q\left(P(\Omega_0, B(x,t')\setminus \overline{B(x,t)}) + Cr^{N-1}h(x,r)\right).
\end{align*}
Letting $t' \rightarrow t$, we get 

\[P(\Omega_0, \overline{B(x,t)}) \leqslant QCr^{N-1}h(x,r).\]
Then, we apply the relative isoperimetric inequality (see \cite[Prop 12.37]{maggi_sets_2012}), which yields
\begin{align*}
\left(\frac{r}{2}\right)^Nh(x,r/2) & \leqslant |\Omega_0 \cap B(x,r/2)| \leqslant |\Omega_0 \cap B(x,t)| \\
&\leqslant \gamma P(\Omega_0, B(x,t))^{\frac{N}{N-1}} \leqslant \gamma (Q C)^{\frac{N}{N-1}} r^{N}h(x,r)^{\frac{N}{N-1}}.
\end{align*}

Otherwise, if $h(x,r) = r^{-N}|B(x,r)\setminus \Omega_0|$, we replace $\Omega_0$ by its complementary set and we obtain the same inequality. In conclusion, there exists some constant $C'>0$, such that
\begin{equation}\label{lem 4.10 id3}
h\left(x,\frac{r}{2}\right) \leqslant C' h(x,r)^{\frac{N}{N-1}}.
\end{equation}

\textit{Step 2.} Following the argument of \cite[Lemma 3.30]{david1998quasiminimal}, we can now prove that if there exist $x \in \partial \Omega_0\cap \mathcal{C}$ and $0<r$ such that $B(x,r)\subset \mathcal{C} $ and $h(x,r) \leqslant \varepsilon$, for some absolute constant $\varepsilon >0$, then either $|\Omega_0 \cap B(x,r/2)|=0$ or $|B(x,r/2)\setminus \Omega_0|=0$. We detail this argument in Lemma~\ref{lem : 3.30}. This contradicts the fact that $x \in \spt \mu_{\Omega_0}$.
Hence, there exists an $\varepsilon>0$ such that for all $x \in \spt \mu_{\Omega_0} \cap \mathcal{C}$ and all $0<r<r_0$, $h(x,r) \geqslant \varepsilon$. In other words 
\[\spt \mu_{\Omega_0}\cap \mathcal{C} \subset  \{ x\in \R^{N} : h(x,r) \geqslant \varepsilon, \forall r>0 \text{ such that } B(x,r) \subset \mathcal{C}\}.\]
 By definition of $\spt \mu_{\Omega_0}$, the converse inclusion is also satisfied, and we have therefore the desired result. 
\end{proof}

\begin{lemme}
\label{lem : 3.30}

There exists a small constant $\varepsilon > 0$, such that for all $x \in \spt \mu_{\Omega_0} \cap \mathcal{C}$ and all $r>0$ such that $B(x,r) \subset \mathcal{C}$, $h(x,r) \geqslant \varepsilon$.  

\end{lemme}

\begin{proof}

As mentioned earlier this lemma is based on \cite[Lemma 3.30]{david1998quasiminimal} (where the study is made on cubes instead of balls). We reason by contradiction. Let $\varepsilon >0$ be a constant, to be fixed later, assume that there exist $x \in \spt \mu_{\Omega_0} \cap \mathcal{C}$ and $r>0$ such that $B(x,r) \subset \mathcal{C}$ satisfying $h(x,r) \leqslant \varepsilon$. The objective is now to prove that

\[\text{ either }|B(x,\frac{r}{2})\cap \Omega_0| = 0 \text{ or } |B(x,\frac{r}{2})\setminus \Omega_0| = 0.\]
By Lebesgue's density theorem, it is enough to show that either all $y\in B(x,\frac{r}{2})$ have 0 density in $\Omega_0$ or they all have 0 density in its the complementary set.

Let $y\in B(x,\frac{r}{2})$. Then, for all $j\geqslant 1 $, $B(y,\frac{r}{2^j}) \subset B(x,r)$. This yields, in particular, that 
\begin{equation}
\label{lem 3.30 initialisation}
h(y,\frac{r}{2}) \leqslant 2^{N} h(x,r) \leqslant 2^{N}\varepsilon.
\end{equation}
\textit{Step 1.} Let us detail why $h(y,\frac{r}{2^j}) \xrightarrow[j\to + \infty]{} 0$. It is necessary to distinguish between two distinct situations. First, if $y \notin \spt \mu_{\Omega_0}$, then, for all $j\geqslant 1$, by \cite[Proposition 12.19]{maggi_sets_2012}, either $\frac{|\Omega_0 \cap B(y,\frac{r}{2^j})|}{|B(y,\frac{r}{2^j})|} = 0$ or $ \frac{|B(y,\frac{r}{2^j}) \setminus \Omega_0|}{|B(y,\frac{r}{2^j})|}=0$. Thus, for all $j \geqslant 1$ $h(y, \frac{r}{2^j})=0$ and the desired limit is obtained. So there only remains to see what happens when $y\in \spt \mu_{\Omega_0}$. In that case, we have seen in \textit{Step 1} of Lemma~\ref{lem 4.10} the inequality \eqref{lem 4.10 id3}, which is valid for any point of the support $\spt \mu_{\Omega_0}$ and any radius $\rho>0$ such that $B(y,\rho) \subset \mathcal{C}$. Hence, for all $j \geqslant 1 $,  
\[h\left(y,\frac{r}{2^{j+1}}\right) \leqslant C' h\left(y,\frac{r}{2^{j}}\right)^{\frac{N}{N-1}}.\]
Thus, with \eqref{lem 3.30 initialisation}, for $\varepsilon$ small enough such that $C'(2^{N} \varepsilon)^{\frac{1}{N-1}} < \frac{1}{2}$, we get 
\[h\left(y,\frac{r}{4}\right)\leqslant C'h\left(y,\frac{r}{2}\right)^{\frac{1}{N-1}}h\left(y,\frac{r}{2}\right)\leqslant C'(2^{N} \varepsilon)^{\frac{1}{N-1}}h\left(y,\frac{r}{2}\right) \leqslant \frac{1}{2}h\left(y,\frac{r}{2}\right) \leqslant 2^{N-1}\varepsilon.\]
So in turn, 
\[h\left(y,\frac{r}{2^{3}}\right)\leqslant C'h\left(y,\frac{r}{4}\right)^{\frac{1}{N-1}}h\left(y,\frac{r}{4}\right)\leqslant C'(2^{N} \varepsilon)^{\frac{1}{N-1}}h\left(y,\frac{r}{4}\right) \leqslant \frac{1}{2}h\left(y,\frac{r}{4}\right) \leqslant 2^{N-2}\varepsilon.\]
And by induction, we get for all $j \geqslant 1$,
\begin{equation}
\label{lem 3.30 4.7}
h\left(y,\frac{r}{2^{j}}\right) \leqslant 2^{N-j+1}\varepsilon,
\end{equation}
which yields the announced limit in this case as well.

\textit{Step 2.} Now, suppose that $h(x,r) = r^{-N} |B(x,r) \cap \Omega_0|$. Then, 
\[\left|B\left(y,\frac{r}{2}\right) \cap \Omega_0 \right | \leqslant|B(x,r) \cap \Omega_0| = r^{N}h(x,r)\leqslant r^{N} \varepsilon = \frac{2^N}{\omega_{N}}\left |B\left(y,\frac{r}{2}\right)\right|\varepsilon. \]
So, for $\varepsilon$ small enough, more precisely if $\frac{2^N}{\omega_{N}}\varepsilon < \frac{1}{2}$, by decomposing the volume of the ball $B(y,\frac{r}{2})$ between the volume of $B(y,\frac{r}{2}) \cap \Omega_0$ and the volume of $B(y,\frac{r}{2}) \setminus \Omega_0$, the previous inequality yields   
\begin{equation}
h\left(y,\frac{r}{2}\right) = \frac{\left|B\left(y,\frac{r}{2}\right) \cap \Omega_0 \right |}{\left|B\left(y,\frac{r}{2}\right) \right |}\ \omega_{N}.
\end{equation}
The same argument holds for $h(y, \frac{r}{4})$ if we assume $\frac{4^N}{\omega_{N}}\varepsilon < \frac{1}{2}$ :
\begin{align*}
\left|B\left(y,\frac{r}{4}\right) \cap \Omega_0 \right | &\leqslant \left |B\left(y,\frac{r}{2}\right) \cap \Omega_0 \right| = \left(\frac{r}{2}\right)^{N}h\left(y,\frac{r}{2}\right)
\leqslant \left( \frac{r}{2} \right)^{N} 2^N \varepsilon = r^N \varepsilon \\
&= \frac{4^N}{\omega_{N}}\left |B\left(y,\frac{r}{4}\right)\right|\varepsilon. 
\end{align*}
Now, we iterate 
\begin{align*}
\left|B\left(y,\frac{r}{8}\right) \cap \Omega_0 \right | &\leqslant \left |B\left(y,\frac{r}{4}\right) \cap \Omega_0 \right| = \left(\frac{r}{4}\right)^{N}h\left(y,\frac{r}{4}\right)\leqslant \left( \frac{r}{4} \right)^{N} 2^{N-1} \varepsilon \\
&= \frac{2^N2^{N-1}}{\omega_{N}}\left |B\left(y,\frac{r}{8}\right)\right|\varepsilon \leqslant \frac{4^N}{\omega_{N}}\left |B\left(y,\frac{r}{8}\right)\right|\varepsilon. 
\end{align*}
Hence, \[h\left(y,\frac{r}{8}\right) = \frac{\left|B\left(y,\frac{r}{8}\right) \cap \Omega_0 \right |}{\left|B\left(y,\frac{r}{8}\right) \right |}\ \omega_{N}.\]
This argument allows us to conclude that for all $j\geqslant 1$, if $\varepsilon$ is chosen small enough so that $\frac{4^N\varepsilon}{\omega_N}<\frac{1}{2}$, then
\begin{equation*}
h\left(y,\frac{r}{2^{j}}\right) = \frac{\left|B\left(y,\frac{r}{2^{j}}\right) \cap \Omega_0 \right |}{\left|B\left(y,\frac{r}{2^{j}}\right) \right |}\ \omega_{N}.
\end{equation*}
Hence, using \eqref{lem 3.30 4.7}, this yields, 
\[\frac{|B(y,\frac{r}{2^{j}})\cap \Omega_0|}{|B(y,\frac{r}{2^{j}})|} \xrightarrow[j \to \infty]{} 0.\]
Thus, all $y~ \in~ B(x,\frac{r}{2})$ have density 0 in $\Omega_0$.

Similarly, if $h(x,r) = r^{-N} |B(x,r) \setminus \Omega_0|$ the same argument shows that all $y \in B(x,r/2)$ have density 0 in $(\Omega_0)^{c}$. This concludes the proof of this Lemma. 
\end{proof}

\begin{lemme}
\label{lem : cond B int}
Let $\Omega_0$ be a Plateau-quasi-minimizer such that $\spt \mu_{\Omega_0} = \partial \Omega_0$, then the open sets 
\begin{align*}
\mathcal{O}_0 &= \{x \in \R^{N} : \exists r> 0 \ |B(x,r)\setminus \Omega_0| = 0 \}\\
\mathcal{O}_1 &= \{x \in \R^{N} : \exists r> 0 \ |B(x,r)\cap \Omega_0| = 0 \}
\end{align*}
 satisfy Condition B. Moreover $\mathcal{O}_0$ is equivalent to $\Omega_0$ and $\spt \mu_{\mathcal{O}_0} = \partial \mathcal{O}_0 = \partial \mathcal{O}_1 = \spt \mu_{\mathcal{O}_1} = \spt \mu_{\Omega_0} \subset \partial \Omega_0$.
\end{lemme}

\begin{proof}
See, for instance, \cite[Lemma 3.6]{rigot2000uniform} for the proof of this lemma. 

Notice that, for any set of finite perimeter $\Omega$ such that $\spt \mu_{\Omega} = \partial \Omega$ with $\partial \Omega$ Ahlfors regular, $\mathcal{O}_0$ (as defined in the Lemma) is equivalent to $\Omega$ and $\spt \mu_{\mathcal{O}_0} = \partial \mathcal{O}_0 = \partial \mathcal{O}_1 = \spt \mu_{\mathcal{O}_1} = \spt \mu_{\Omega} \subset \partial \Omega$. Therefore, the Plateau-quasi-minimality assumption is only necessary to justify that the boundary is Ahlfors regular and to prove that $\mathcal{O}_0$ satisfies Condition B.
\end{proof}

\begin{proof}{\textit{(Proposition~\ref{prop : rect int})}}
The previous Lemmas lead to the existence of an open set $\mathcal{O}_0$ in the equivalence class of $\Omega_0$, which satisfies Condition B inside the convex set $\mathcal{C}$. This achieves the proof of Proposition~\ref{prop : rect int}.
\end{proof}


\subsection{Uniform rectifiability at the boundary}

Now, we proceed to extend this result to the boundary. As explained before, thanks to Proposition \ref{prop : caract unif rect}, it only remains to show that Condition B is also satisfied at the boundary. 
\begin{prop}[Condition B at the boundary]
\label{prop : rect bdd}
 Let $(\mathcal{D}, \mathcal{C}, E_0)$ satisfying Hypothesis H (see Definition~\ref{def H}) and let $\Omega_{0}$ be a Plateau-quasi-minimizer. Then, there exist an equivalent open set with the same assumption, still denoted by $\Omega_0$, and a constant $C>1$ such that for all $x\in \partial \Omega_{0}\cap \mathcal{D}$ and $r>0$  satisfying $B(x,r)  \subset \mathcal{D}$ and $B(x,r) \cap \partial \mathcal{C} \neq \emptyset$, there exist two balls $B_1 \subset B(x,r)\cap \Omega_0$ and $B_2 \subset B(x,r) \setminus \overline{\Omega_0}$ with radius larger than $C^{-1}r$. 
\end{prop}

\begin{proof}

We have established that $\mathcal{O}_0$, defined in Lemma~\ref{lem : cond B int}, is an open set equivalent to $\Omega_0$, such that $\spt \mu_{\mathcal{O}_0} = \partial \mathcal{O}_0$. Thus, since $\Omega_0$ is a Plateau-quasi-minimizer, by Remark~\ref{rq : equivalent sets}, $\mathcal{O}_0$ is also a Plateau-quasi-minimizer. Hence, by the previous section, its boundary $\partial \mathcal{O}_0$ is Ahlfors regular.Therefore, there remains to show that $\mathcal{O}_0$ satisfies Condition B. 

As it was done for the Ahlfors regularity at the boundary we distinguish differents situations, in which we apply either Condition B' or the interior regularity argument. It is important to note that one of the key arguments, which ensures the extension of the uniform rectifiability result to the boundary, is the inequality for P-quasi-minimizers, which do not cross the surface $\Sigma$, proven for the Ahlfors regularity result in Lemma~\ref{lem AR},
\[P(\mathcal{O}_0,B(x,t')) \leqslant Q P(\mathcal{O}_0\setminus \overline{B(x,t)},B(x,t')).\]

\begin{enumerate}
\item First, if $d(x,\Sigma) \leqslant \frac{r}{2}$ we use Condition B'. Indeed, there exists $y\in \Sigma \cap B(x,\frac r 2)$ and Condition B' yields the existence of two balls $B_1 \subset B(y,r/2)\cap (E_0 \setminus \overline{\mathcal{C}})$ and $B_2 \subset (B(y,r/2)\setminus \overline{\mathcal{C}}) \setminus \overline{E_0} $ with radii larger than $C^{-1} r/2$. Thus, by boundary constraint it follows that $B_1 \subset B(x,r)\cap \mathcal{O}_0$ and $B_2 \subset  B(x,r) \setminus \overline{\mathcal{O}_0}$. This achieves the proof of the Proposition for this case.

\item Otherwise, if $d(x,\Sigma) > \frac{r}{2}$, we use the interior result in the ball $B(x,\frac{r}{2})$. In this case, the sole difference with the proof of the interior rectifiability is in Lemma~\ref{lem 4.10}. In this proof we use that since $B(x,r) \subset \mathcal{C}$ the set $\mathcal{O}_0\setminus \overline{B(x,t)}$ (with $t<r$) is a competitor. This is no longer systematically the case at the boundary. Nevertheless, if $\mathcal{O}_0\setminus \overline{B(x,t)}$ (with $t<\frac{r}{2}$) is not a competitor, it was proven in Lemma~\ref{lem AR} that the inequality
\[P(\mathcal{O}_0,B(x,t')) \leqslant QP(\mathcal{O}_0\setminus \overline{B(x,t)},B(x,t'))\] (with $t<t'<\frac{r}{2}$) remains satisfied in this case. Given this information, the remainder of the proof, previously developed for the interior case, can be carried out in the same manner. 
\end{enumerate}

\end{proof}

\begin{proof}{\textit{(Theorem~\ref{th : unif rect})}}
Propositions~\ref{prop : rect int} and \ref{prop : rect bdd} insure that Condition B is satisfied up to the boundary. Moreover, we have proven the existence of an open set $\mathcal{O}_0$ in the equivalence class of $\Omega_0$, which satisfies Condition B inside $\mathcal{D}$ and has a regular boundary, in the sense of Definition~\ref{def : regbdd}.
 Thus, we apply Proposition~\ref{th 6.1 David} to achieve the proof of Theorem~\ref{th : unif rect}.
\end{proof}

\section{Characterization of Plateau-quasi-minimizer by bi-John domains}
\label{section John}
The goal of this section is to show that every Plateau-quasi-minimizer is a bi-John domain with regular boundary, in the sense of Definition~\ref{def : regbdd}, and that the converse is also true. In other words, we can characterize Plateau-quasi-minimizer by a purely geometric property.

In \cite{david1998quasiminimal} David and Semmes used John domains to characterize quasiminimal sets. In this section, we will adapt their result to Plateau's problem. Namely, we prove that their characterization remains valid if we add a boundary condition. For the reader's convenience, we recall the statement of our main Theorem, presented in the Introduction.

\begin{thm}[Characterization of Plateau-quasi-minimizer by bi-John domains]
\label{th john dom}
Let $(\mathcal{D}, \mathcal{C}, E_0)$ satisfying Hypothesis H (see Definition~\ref{def H}) and let $\Omega_0$ be a competitor (see Definition~\ref{competitors}).
Then $\Omega_0$ is a Plateau-quasi-minimizer if and only if there exists an equivalent open set, which is a bi-John domain, with regular boundary in the sense of Definition~\ref{def : regbdd}.
\end{thm}

\subsection{The direct implication}
\label{section direct implication}
We start by proving the direct implication of Theorem~\ref{th john dom}, \emph{i.e.}, we show that a Plateau-quasi-minimizer is a bi-John domain. To that aim, we introduce the notion of domain of isoperimetry, coming from \cite{david1998quasiminimal}. 

\begin{definition}[Domain of isoperimetry]
\label{def : domain isop}
An open set $\Omega \subset \R^{N}$ is a domain of isoperimetry if there exists a constant $C>0$ such that for all open set $\omega \subset \Omega$, 
\begin{equation}
\min(|\omega|,|\Omega \setminus \omega|) \leqslant C \Haus(\partial \omega \cap \Omega)^{\frac{N}{N-1}}.
\end{equation}
\end{definition}

\begin{prop}{\cite[Theorem 6.1]{david1998quasiminimal}}
\label{th 6.1 David}
Let $\Omega$ be a bounded domain that contains $B(0,1)$. Assume that $\Omega$ satisfies Condition B (see Definition~\ref{condition B}) and is a domain of isoperimetry. Then $\Omega$ is a John domain with center $0$. The constant for the John condition is controlled by $diam(\Omega)$ and the constants from Condition B and domain of isoperimetry. 
\end{prop}

Thus, if we can prove that a Plateau-quasi-minimizer, satisfying Condition B, is a domain of isoperimetry then, Proposition~\ref{th 6.1 David} yields the direct implication of Theorem~\ref{th john dom}. In Section \ref{section unif rect}, we have established that for any P-quasi-minimizer there exists a open equivalent set $\mathcal{O}_0$ such that $\spt \mu_{\Omega_0} = \partial \Omega_0 $ that satisfies Condition B up to the boundary. Therefore, it is enough to show that such a Plateau-quasi-minimizer is a domain of isoperimetry. To do so, we need to characterize domains of isoperimetry. This is done in Lemma~\ref{lem : domain isop}.

\begin{lemme}
\label{lem : domain isop}
Let $\Omega \subset \R^{N}$ be an open set such that the boundary condition, $\Omega \setminus \overline{\mathcal{C}} = E_0 \setminus \overline{\mathcal{C}}$  is satisfied. If there exists a constant $C>0$ such that for all open set $\tilde{\omega} \subset \Omega$ satisfying the boundary condition ($\tilde \omega \setminus \overline{\mathcal{C}} = E_0 \setminus \overline{\mathcal{C}}$), 
\begin{equation}
\label{eq : domain isop}
|\Omega \setminus \tilde{\omega}| \leqslant C P(\tilde{\omega}, \Omega)^{\frac{N}{N-1}},
\end{equation}
then $\Omega$ is a domain of isoperimetry.
\end{lemme}

\begin{proof}

Assume that $\Omega$ satisfies the boundary condition and \eqref{eq : domain isop} for all open sets $\tilde{\omega} \subset \Omega$ such that $\tilde \omega \setminus \overline{\mathcal{C}} = E_0 \setminus \overline{\mathcal{C}} =: \mathcal{A}^+$. Let $\omega \subset \Omega$ be an open set. We assume that $\omega$ does not satisfy the boundary condition, \emph{i.e.}, $\Omega \setminus \overline{\mathcal{C}} \not\subset \omega \setminus \overline{\mathcal{C}}$, otherwise the result is trivial. 

We distinguish two possible cases. First, if 
\begin{equation}
\label{lem 5.5 cas 1}
|\mathcal{A}^+ \setminus \omega | = \min(|\omega |,|\mathcal{A}^+ \setminus \omega|).
\end{equation}
We introduce the open set $\tilde{\omega} = \omega \cup \mathcal{A}^+ \subset \Omega$, which satisfies the boundary condition. Notice, that we then have the following decomposition 
\[|\Omega \setminus \omega| = |\Omega \setminus \tilde{\omega}|+|\mathcal{A}^+ \setminus \omega|.\]
Then, because $\mathcal{A}^+$ is a domain of isoperimetry, thanks to Hypothesis~\ref{H1}, it follows from the assumption \eqref{lem 5.5 cas 1} that
\begin{align*}
|\Omega \setminus \omega| &\leqslant |\Omega \setminus \tilde{\omega}| + C\Haus(\mathcal{A}^+\cap \partial \omega)^{\frac{N}{N-1}} \\
&\leqslant |\Omega \setminus \tilde{\omega}| + C\Haus(\Omega \cap \partial \omega)^{\frac{N}{N-1}},
\end{align*}
since $\mathcal{A}^+ \subset \Omega$. 
Moreover, $\tilde{\omega}$ satisfies the boundary condition. Thus,
\[|\Omega \setminus \tilde{\omega}| \leqslant C P(\tilde{\omega},\Omega)^{\frac{N}{N-1}}.\]
Therefore, there only remains to show that 
\[P(\tilde{\omega}, \Omega)\leqslant C \Haus(\partial \omega \cap \Omega).\]
Furthermore, Lemma \ref{Maggi_cvx} yields
\begin{align*}
\Haus(\partial \omega \cap \Omega) \geqslant \Haus(\partial^* \omega \cap \Omega) = P(\omega,\Omega) = P(\Omega \setminus \omega, \Omega) 
 \geqslant P(\Omega \setminus \omega \cap \overline{\mathcal{C}}, \Omega),
\end{align*}
because $\partial \omega \supset \partial^* \omega$.
Notice that since $\omega \cap \overline{\mathcal{C}} = \tilde\omega \cap \overline{\mathcal{C}}$ and $ \Omega\setminus \tilde{\omega} \subset \overline{\mathcal{C}}$ we get 
\begin{align*}
\Haus(\partial \omega \cap \Omega)
 \geqslant P(\Omega \setminus \omega \cap \overline{\mathcal{C}}, \Omega) = P(\Omega \setminus \tilde{\omega}, \Omega)= P(\tilde{\omega},\Omega).
\end{align*}
This achieves the proof in this case.

Otherwise, 

\begin{equation}
\label{lem 5.5 cas 2}
|\omega | = \min(|\omega |,|\mathcal{A}^+ \setminus \omega|).
\end{equation}
In that case, we argue with $\omega' = \Omega \setminus \omega$ as previously. Indeed, 
\[|\mathcal{A}^+ \setminus \omega'| = |\mathcal{A}^+ \cap \omega| \leqslant |\omega| \leqslant |\mathcal{A}^+ \setminus \omega| \leqslant |\Omega \setminus \omega| = |\omega'|.\]
Thus, the previous case leads to 
\[|\omega| = |\Omega \setminus \omega'| \leqslant C \Haus(\partial \omega' \cap \Omega)^{\frac{N}{N-1}} = C \Haus(\partial \omega \cap \Omega)^{\frac{N}{N-1}},\]
which achieves the proof for this second case.
\end{proof}

With the previous Lemma, we can prove the desired result. 

\begin{prop}
\label{prop : dom isop}
Let $(\mathcal{D}, \mathcal{C}, E_0)$ satisfying Hypothesis H (see Definition~\ref{def H}) and let $\Omega$ be an open Plateau-quasi-minimizer such that $\spt \mu_{\Omega} = \partial \Omega$. Then $\Omega$ is a domain of isoperimetry.
\end{prop}

\begin{proof}
To prove that $\Omega$ is a domain of isoperimetry, we use the caracterization of Lemma~\ref{lem : domain isop}. In other words, it is enough to show that for all $\omega \subset \Omega$ open and satisfying the boundary condition we have 
\[|\Omega \setminus \omega| \leqslant C P(\omega, \Omega)^{\frac{N}{N-1}}.\]
We can assume that $P( \omega, \Omega)^{\frac{N}{N-1}}$ is finite, otherwise the inequality is trivial. In particular, this implies that $\chi_\omega \in BV(\mathcal{D})$, hence $\omega$ is a competitor for Plateau's problem. 

Besides, the isoperimetric inequality (\cite[Theorem 14.1]{maggi_sets_2012}) ensures the existence of a constant $C>0$, depending only on the dimension, such that 
\[|\Omega \setminus \omega| \leqslant C P(\Omega \setminus \omega)^{\frac{N}{N-1}}.\]
Therefore it is enough to show that 
\begin{equation}
P(\Omega\setminus \omega) \leqslant C' P( \omega, \Omega).
\end{equation}
Let us prove this inequality. First, since $\omega \subset \Omega$, \eqref{Maggi prive} yields
\begin{align*}
P(\Omega \setminus \omega) &= P(\Omega, \omega^{(0)}) + P(\omega, \Omega^{(1)}) + \Haus(\{\nu_\Omega = - \nu_\omega\}) \\
&= P(\Omega, \omega^{(0)}) + P(\omega, \Omega^{(1)})
\end{align*}
Furthermore, since $\Omega$ is a P-quasi-minimizer such that $\spt \mu_{\Omega} = \partial \Omega$, Proposition~\ref{th : AhlReg bdd} yields that $\partial \Omega = \partial^* \Omega $ $\Haus$-a.e. and Proposition~\ref{prop : rect bdd} implies in particular that $\Omega$ is open (up to considering an equivalent set). Thus, $\Omega^{(1)} = \Omega$ $\Haus$-a.e., and it remains only to show 
\[P(\Omega, \omega^{(0)}) \leqslant C''P(\omega, \Omega).\]
We use the decomposition, $\omega^{(0)} = \R^{N} \setminus (\omega^{(1)} \cup \partial^*\omega)$. 
\begin{align*}
P(\Omega, \omega^{(0)}) &= \Haus(\partial^*\Omega \cap \omega^{(0)})\\
&= \Haus(\partial^*\Omega \setminus (\omega^{(1)} \cup \partial^*\omega))\\
&\leqslant \Haus(\partial^*\Omega \setminus\partial^*\omega).
\end{align*}
Since $\Omega$ is a P-quasi-minimizer and $\omega$ a competitor it follows
\[P(\Omega, \omega^{(0)}) \leqslant Q\Haus(\partial^*\omega \setminus\partial^*\Omega).\]
Finally, since $\omega \subset  \Omega$, the same decomposition as earlier  implies that $\partial^*\omega \cap \Omega^{(0)} = \emptyset$,
\[P(\Omega, \omega^{(0)}) \leqslant Q\Haus(\partial^*\omega \cap\Omega^{(1)}) = QP(\omega, \Omega).\]
 This yields 
\begin{align*}
P(\Omega \setminus \omega) \leqslant (Q+1)P(\omega, \Omega). 
\end{align*}
We conclude with Lemma~\ref{lem : domain isop} that $\Omega_0$ is a domain of isoperimetry.
\end{proof}

We follow the same strategy to show that $\mathcal{D}\setminus \overline{\Omega}$ is also a domain of isoperimetry.

\begin{prop}
\label{prop : dom isop2}
Let $(\mathcal{D}, \mathcal{C}, E_0)$ satisfying Hypothesis H (see Definition~\ref{def H}) and let $\Omega$ be an open Plateau-quasi-minimizer such that $\spt \mu_{\Omega} = \partial \Omega$. Then $\mathcal{D}\setminus \overline{\Omega}$ is a domain of isoperimetry.
\end{prop}

\begin{proof}
By symmetry of the problem we could consider $\mathcal{D}\setminus \overline{E_0}$ as boundary constraint instead of $E_0$. And, if $\Omega$ is an open P-quasi-minimizer such that $\spt \mu_{\Omega_0} = \partial \Omega_0$ for the problem with boundary constraint $E_0$ then $ \mathcal{D}\setminus \overline{\Omega}$ is an open P-quasi-minimizer of Plateau's problem with boundary constraint $\mathcal{D}\setminus \overline{E_0}$,  such that $\spt \mu_{\mathcal{D}\setminus \overline{\Omega}} = \partial (\mathcal{D}\setminus \overline{\Omega})$. Therefore the same proof unfolds since the same assumptions are made for $\mathcal{D}\setminus \overline{E_0}$ as for $E_0$.
\end{proof}

\begin{proof}{(\emph{Direct implication of Theorem~\ref{th john dom}})}
 Theorem~\ref{th : AhlReg} and Theorem~\ref{th : unif rect} yield that, up to considering a proper equivalent set, $\Omega$ satisfies Condition B (see Definition~\ref{condition B}). And, Proposition~\ref{prop : dom isop} implies that a $\Omega$ is a domain of isoperimetry. Hence, by Proposition~\ref{th 6.1 David} $\Omega$ is a John domain in $\mathcal{D}$. 
 
 Indeed, the assumption $B(0,1) \subset \Omega_0$ in Proposition \ref{th 6.1 David} can be easily replaced by the assumption $B(x_0,r) \subset \Omega_0$ where $x_0 \in \mathcal{A}^{+}$ and $r>0$ such that $B(x_0,r) \subset \mathcal{A}^+$, which is always satisfied for a P-quasi-minimizer, since it is a competitor, and thus contains $\mathcal{A}^+$. 
 
 Similarly, Proposition~\ref{prop : dom isop2} combined with Proposition~\ref{th 6.1 David} yield that there exists an open set equivalent to the complementary $\mathcal{D}\setminus \overline{\Omega}$, which is a John domain in $\mathcal{D}$. This leads to the direct implication of Theorem~\ref{th john dom}. 
\end{proof}

\subsection{The converse implication }

 Let $\Omega_0$ be a competitor which is a bi-John domain with regular boundary, in the sense of Definition~\ref{def : regbdd}. In particular, this yields that $\spt \mu_{\Omega_0} = \partial\Omega_0$. We denote $z_0$ the center of $\Omega_0$ and $C>0$ the John-constant.

Inspired by \cite{david1998quasiminimal}, we introduce a weighted Plateau problem, and we show its well-posedness. 

\begin{definition}

Let $A\geqslant 1$ be a constant, and $\Omega_0$ be as defined in the beginning of the Section, then we define the associated weighted Plateau problem as 
\begin{equation}
\label{def Pa}
\min \{\Haus(\partial^* \Omega \cap \partial \Omega_0\cap \mathcal{D})+A\Haus(\partial^*\Omega \setminus \partial \Omega_0\cap \mathcal{D}), \text{ for } \Omega \text{ competitor } \}.
\end{equation}

\end{definition}

\begin{prop}
The weighted Plateau problem \eqref{def Pa} admits solutions.
\end{prop}

\begin{proof}
The proof of this Proposition is similar to the proof of Proposition~\ref{existence} (the existence of solutions of Plateau's problem). It uses the direct method of Calculus of variations, with the compactness result for set of finite perimeter and the lower semi continuity of the perimeter in an open set (see for instance \cite[Proposition 12.15]{maggi_sets_2012}). However, we need to explain how we pass the first term $\Haus(\partial^*\Omega_n\cap \partial \Omega_0\cap \mathcal{D})$ to the liminf, as we look at a perimeter in a closed set. 

Let $(\Omega_n)$ be a minimizing sequence for \eqref{def Pa}, which converges, by Rellich Theorem, up to extracting a subsequence, to $\Omega$ a competitor. We will justify why 
\begin{align*}
&\liminf (\Haus(\partial^* \Omega_n\cap \partial\Omega_0\cap \mathcal{D})+A\Haus(\partial^* \Omega_n\setminus \partial\Omega_0\cap \mathcal{D}) )\\
&\geqslant \Haus(\partial^* \Omega\cap \partial\Omega_0\cap \mathcal{D})+A\Haus(\partial^* \Omega\setminus \partial\Omega_0\cap \mathcal{D}).
\end{align*}

Assume that there exists an open set $U\subset \mathcal{D}$ such that $\partial \Omega_0 \subset U$, $\Haus(\partial^* \Omega_n\cap \partial U)=0$ and $\Haus(\partial^* \Omega \cap \partial U)=0$. Then, we decompose
\begin{align*}
&\Haus(\partial^* \Omega_n\cap \partial\Omega_0\cap \mathcal{D})+A\Haus(\partial^* \Omega_n\setminus \partial\Omega_0\cap \mathcal{D}) \\
&= \Haus(\partial^* \Omega_n\cap U\cap \mathcal{D} )+A\Haus(\partial^* \Omega_n \cap \mathcal{D}\setminus\overline{U}) + (A-1)\Haus(\partial^* \Omega_n\cap U \setminus \partial \Omega_0).
\end{align*} 
Since $U$ , $\mathcal{D}\setminus\overline{U}$ and $U \setminus \partial \Omega_0$ are open sets, the lower semi continuity of the perimeter yields 
\begin{align*}
& \liminf (\Haus(\partial^* \Omega_n\cap \partial\Omega_0\cap \mathcal{D})+A\Haus(\partial^* \Omega_n\setminus \partial\Omega_0\cap \mathcal{D})) \\
&= \liminf( \Haus(\partial^* \Omega_n\cap U\cap \mathcal{D} )+A\Haus(\partial^* \Omega_n \cap \mathcal{D}\setminus\overline{U})\\
& \qquad \qquad \qquad \qquad  \qquad \qquad \qquad+ (A-1)\Haus(\partial^* \Omega_n\cap U \setminus \partial \Omega_0))\\
&\geqslant \liminf( \Haus(\partial^* \Omega_n\cap U\cap \mathcal{D} ))+A\liminf(\Haus(\partial^* \Omega_n \cap \mathcal{D}\setminus\overline{U})) \\
& \qquad \qquad \qquad \qquad  \qquad \qquad \qquad + (A-1)\liminf(\Haus(\partial^* \Omega_n\cap U \setminus \partial \Omega_0)))\\
&\geqslant \Haus(\partial^* \Omega \cap U\cap \mathcal{D} )+A\Haus(\partial^* \Omega \cap \mathcal{D}\setminus\overline{U}) + (A-1)\Haus(\partial^* \Omega\cap U \setminus \partial \Omega_0)\\
&= \Haus(\partial^* \Omega\cap \partial\Omega_0\cap \mathcal{D})+A\Haus(\partial^* \Omega\setminus \partial\Omega_0\cap \mathcal{D}).
\end{align*}
To conclude, we justify the existence of such a $U$. The coarea formula applied with the 1-Lipschitz function $f = \mathrm{dist(\ \cdot \ ,\partial \Omega_0)}$ to the rectifiable set $\partial^*\Omega_n\cap\mathcal{D}\cap \{f < 1\}$ implies that the integral $\int_0^1 \mathcal{H}^{N-2}(\partial^*\Omega_n \cap\mathcal{D}\cap \{f=t\})dt $ is finite. Thus, for almost every $t \in [0,1]$, $\mathcal{H}^{N-2}(\partial^*\Omega_n \cap\mathcal{D}\cap \{f=t\}) < \infty$, which leads to $\Haus(\partial^*\Omega_n \cap\mathcal{D}\cap \{f=t\}) = 0$. Hence, we can consider $U = \{f < t\}$, with $t$ satisfying the condition above for all $\Omega_n$ and for $\Omega$, which is possible since the family $(\Omega_n)$ is countable.  
\end{proof}

Now, we explain the connection between Plateau's problem and the weighted Plateau problem \eqref{def Pa}. 

\begin{prop}
\label{prop : min implique quasi}
Let $\Omega$ be a minimizer of the weighted Plateau problem \eqref{def Pa}.
Then, $\Omega$ is a Plateau-quasi-minimizer.
\end{prop}

\begin{proof}
Let $\Tilde{\Omega}$ be a competitor. We must prove that there exists a constant $Q \geqslant 1$, independent of $\Tilde{\Omega}$, such that 
\[\Haus(\partial^*\Omega \setminus \partial^*\Tilde{\Omega} \cap \mathcal{D}) \leqslant Q \Haus(\partial^*\Tilde{\Omega} \setminus \partial^*\Omega \cap \mathcal{D}).\]

Since $\Omega$ is a solution of \eqref{def Pa} and $\Tilde{\Omega}$ a competitor this yields
\begin{align*}
AP(\Omega,\mathcal{D}) &= A\Haus(\partial^*\Omega \cap \partial \Omega_0 \cap \mathcal{D}) + A\Haus(\partial^*\Omega \setminus \partial \Omega_0\cap \mathcal{D})\\
&= \Haus(\partial^*\Omega \cap \partial \Omega_0\cap \mathcal{D}) + A\Haus(\partial^*\Omega \setminus \partial \Omega_0\cap \mathcal{D})\\
 & \quad + (A-1)\Haus(\partial^*\Omega \cap \partial \Omega_0\cap \mathcal{D})\\
&\leqslant \Haus(\partial^*\Tilde{\Omega} \cap \partial \Omega_0\cap \mathcal{D}) + A\Haus(\partial^*\Tilde{\Omega} \setminus \partial \Omega_0\cap \mathcal{D}) \\
& \quad + (A-1)\Haus(\partial^*\Omega \cap \partial \Omega_0\cap \mathcal{D})\\
&= AP(\Tilde{\Omega},\mathcal{D}) + (A-1)(\Haus(\partial^*\Omega \cap \partial \Omega_0\cap \mathcal{D}) - \Haus(\partial^*\Tilde{\Omega} \cap \partial \Omega_0\cap \mathcal{D})).
\end{align*}
Since $A \geqslant 1$, decomposing the perimeter leads to

\begin{align*}
A\Haus(\partial^*\Omega \setminus \partial^*\Tilde{\Omega} \cap \mathcal{D}) + &A\Haus(\partial^*\Omega \cap \partial^*\Tilde{\Omega} \cap \mathcal{D}) \leqslant A\Haus(\partial^*\Tilde{\Omega} \setminus \partial^*\Omega \cap \mathcal{D})\\
&  +A\Haus(\partial^*\Tilde{\Omega} \cap \partial^*\Omega \cap \mathcal{D})\\
& + (A-1)\Haus(\partial \Omega_0 \cap \partial^*\Omega \setminus \partial^* \Tilde{ \Omega} \cap \mathcal{D}).
\end{align*}
Thus, 
\begin{align*}
A\Haus(\partial^*\Omega \setminus \partial^*\Tilde{\Omega} \cap \mathcal{D}) &\leqslant A\Haus(\partial^*\Tilde{\Omega} \setminus \partial^*\Omega \cap \mathcal{D})\\
& + (A-1)\Haus(\partial \Omega_0 \cap \partial^*\Omega \setminus \partial^* \Tilde{ \Omega} \cap \mathcal{D})\\
&\leqslant A\Haus(\partial^*\Tilde{\Omega} \setminus \partial^*\Omega \cap \mathcal{D}) + (A-1)\Haus(\partial^*\Omega \setminus \partial^* \Tilde{ \Omega} \cap \mathcal{D}).
\end{align*}
Finally, this allows us to conclude the proof of this Proposition
\[\Haus(\partial^*\Omega \setminus \partial^*\Tilde{\Omega} \cap \mathcal{D})  \leqslant A\Haus(\partial^*\Tilde{\Omega} \setminus \partial^*\Omega \cap \mathcal{D}).\]
\end{proof}

With Proposition~\ref{prop : min implique quasi} at hand, finding a constant $A$, such that the bi-John domain considered is solution of the associated weighted Plateau problem, is enough to prove the converse implication of Theorem~\ref{th john dom}. 

\begin{proof}{(\emph{Converse implication of Theorem~\ref{th john dom}})}
Let $\Omega_0$ be a competitor, which is a bi-John domain with regular boundary. We fix a constant $A \geqslant 1$, to be determined later, and $\Omega$ a solution of the associated weighted Plateau problem \eqref{def Pa}. Because $\Omega$ is a minimizer of \eqref{def Pa}, Proposition~\ref{prop : min implique quasi} yields that $\Omega$ is a P-quasi-minimizer. Then, by the direct implication of Theorem~\ref{th john dom} proved in Section~\ref{section direct implication}, $\Omega$ (or an equivalent set) is a bi-John domain with regular boundary. Now, because $\Omega_0$ is a competitor and $\Omega$ is a minimizer of \eqref{def Pa}, it yields
\begin{align*}
\Haus(\partial \Omega \cap \partial \Omega_0 \cap \mathcal{D}) +& A\Haus(\partial\Omega \setminus \partial \Omega_0\cap \mathcal{D}) \leqslant \Haus(\partial \Omega_0 \cap \mathcal{D})\\
&= \Haus(\partial \Omega_0 \cap \partial \Omega \cap \mathcal{D}) + \Haus(\partial \Omega_0 \setminus \partial \Omega \cap \mathcal{D}).
\end{align*}
Hence,
\begin{equation}
\label{ineg avec A}
A\Haus(\partial \Omega \setminus \partial \Omega_0 \cap \mathcal{D}) \leqslant \Haus(\partial\Omega_0 \setminus \partial \Omega \cap \mathcal{D}) .
\end{equation}
We claim that there exist a constant $C>0$, independent of $\Omega$ and $A$ such that 
\begin{align*}
&\Haus(\partial \Omega_0 \setminus \overline{\Omega}\cap \mathcal{D}) \leqslant C \Haus(\partial \Omega \cap \Omega_0\cap \mathcal{D}) \\
&\Haus(\partial \Omega_0 \cap \Omega\cap \mathcal{D}) \leqslant C \Haus(\partial \Omega \setminus \overline{\Omega_0}\cap \mathcal{D}).
\end{align*}
This is proved in Lemma~\ref{lem John reciproque}. Assume for now that this is true. Combining \eqref{ineg avec A} with the claim of Lemma~\ref{lem John reciproque} yields 
\begin{align*}
\Haus(\partial\Omega_0 \setminus \partial \Omega \cap \mathcal{D}) &= \Haus(\partial \Omega_0 \setminus \overline{\Omega}\cap \mathcal{D}) + \Haus(\partial \Omega_0 \cap \Omega\cap \mathcal{D}) \\
&\leqslant C (\Haus(\partial \Omega \cap \Omega_0\cap \mathcal{D}) + \Haus(\partial \Omega \setminus \overline{\Omega_0}\cap \mathcal{D}))\\
&= C \Haus(\partial \Omega  \setminus \partial \Omega_0 \cap \mathcal{D})\\
&\leqslant \frac{C}{A}\Haus(\partial \Omega_0  \setminus \partial \Omega \cap \mathcal{D}).
\end{align*}
We choose $A>C$ so that $\Haus(\partial\Omega_0 \setminus \partial \Omega \cap \mathcal{D})=0$. And with \eqref{ineg avec A} this yields $\Haus(\partial\Omega \setminus \partial \Omega_0 \cap \mathcal{D})=0$. Hence, $\Haus$ almost everywhere in $\mathcal{D}$, $\partial \Omega = \partial \Omega_0$, and since the boundaries are regular it yields $\partial \Omega = \partial \Omega_0$. Then, because $\Omega$ and $\Omega_0$ are bi-John domain it follows that $\R^{N} \setminus \partial \Omega $ and $\R^{N} \setminus \partial \Omega_0$ have exactly two connected components. Finally, since both $\Omega$ and $\Omega_0$ satisfy the boundary condition, we get 
\[\Omega = \Omega_0.\]
Thus, because $\Omega$ is a P-quasi-minimizer, so is $\Omega_0$. 
\end{proof}

There only remains to justify the claim of Lemma~\ref{lem John reciproque}. To that aim, we recall a result from \cite{david1998quasiminimal}, that we will use to prove this result.

\begin{lemme}{\cite[Lemma 7.12]{david1998quasiminimal}}
\label{lem 7.12 DS}
Let $W$ be a bounded John domain with regular boundary, which is a competitor for Plateau's problem, and $f \in BV(W)$.
 We denote $z_0$ its center and $C>0$ the John constant associated. This yields the existence of $r_0 \geqslant C^{-1}$ such that $B_0 := B(z_0,r_0) \subset W$.
 
  For $x \in W$, we denote by $\delta(x)$ the distance between $x$ and the complementary set $\R^N \setminus W$ and we define, for $M \geqslant 1$ a constant depending on $W$,
\[f_{\sharp}(x) = \frac{1}{\delta(x)^N} \int_{B(x,\frac{\delta(x)}{M})}|f(y) - m_0f|dy,\]
where $m_ 0f = \frac{1}{|B_0|}\int_{B_0}f(y)dy$. 
And, for $z \in \partial W$ we define 
\[Nf(z) = \sup \{f_\sharp (x), \ x\in W \text{ such that } |x-z| \leqslant M \delta(x)\}.\]
Then, there exists a constant $C'$, depending only on $W$ such that
\[\int_{\partial W}Nf d\Haus \leqslant C' \int_W |\nabla f|.\]
\end{lemme}

\begin{lemme}
\label{lem John reciproque} 
Let $\Omega$ and $\Omega_0$ be two competitors which are bi-John domains with regular boundary (see Definition~\ref{def : regbdd}). Then there exist a constant $C>0$, only depending on $\Omega_0$, the dimension and the boundary conditions such that 
\begin{align*}
&\Haus((\partial \Omega_0 \setminus \overline{\Omega})\cap \mathcal{D}) \leqslant C \Haus((\partial \Omega \cap \Omega_0)\cap \mathcal{D}) \\
&\Haus((\partial \Omega_0 \cap \Omega)\cap \mathcal{D}) \leqslant C \Haus((\partial \Omega \setminus \overline{\Omega_0})\cap \mathcal{D}).
\end{align*}
\end{lemme}

\begin{proof}
\textit{Step 1. } We start by proving the first inequality. We consider the function $f = \chi_{\Omega_0 \setminus \Omega} \in BV(\Omega_0)$ and we recall the definition introduced in the Lemma~\ref{lem 7.12 DS}. $\Omega_0$ is a John domain with center, denoted by $z_0$, and constant, denoted by $C$. Since $\mathcal{A}^+$ is contained in $\Omega_0$, we can assume, without loss of generality, that $z_0 \in \mathcal{A}^+$. 

Indeed, a John domain with center, as defined in~\ref{def John}, also called a John domain in the carrot sense is a John domain in the cigar sense (this definition involves no center), see \cite[Section 2.17]{vaisala1988uniform} for precise definition of John domain in the carrot and cigar sense. Conversely, a bounded John domain in the cigar sense is a John domain in the carrot sense. This equivalence in shown in \cite[Theorem 2.21]{vaisala1988uniform}. 

Thus, there exists a radius $r_0 \geqslant \Tilde{C}^{-1}$, with $\Tilde{C}$ a constant, such that $B_0 := B(z_0,r_0) \subset \mathcal{A}^+\subset \Omega_0$. Then, we introduce the function $f_\sharp$, defined on $\Omega_0$ and $Nf$ defined on $\partial \Omega_0$ as in Lemma~\ref{lem 7.12 DS}. Here, because $f = \chi_{\Omega_0 \setminus \Omega} $, we have $m_0f = 0$, since the boundary conditions force $B_0 \subset \mathcal{A}^+ \subset \Omega_0$ and $B_0 \subset \mathcal{A}^+ \subset \Omega$. Thus, in this case, since by definition of $\delta(x)$, $B(x,\frac{\delta(x)}{M}) \subset \Omega_0$, it follows 
\[f_{\sharp}(x) =\frac{1}{\delta(x)^N} \left|B\left(x,\frac{\delta(x)}{M}\right) \setminus \Omega\right|.\]
We claim that there exists a constant $C>0$ such that \[\text{for all } z \in \partial \Omega_0 \setminus \overline{\Omega},\ Nf(z) \geqslant C^{-1}.\]
Let us prove this claim. Let $z \in \partial \Omega_0 \setminus \overline{\Omega}$. We denote $\delta_z = dist(z,\overline{\Omega}) >0$. Then, thanks to the John condition, we can find a point $x \in B(z, M\delta(x))\cap B(z,\frac{\delta_z}{2})\cap \Omega_0$. Indeed, the John condition guaranties the existence of a $C$-Lipschitz curve, denoted by $\alpha$, connecting $z$ to the center of $\Omega_0$. So we will choose $x$ in the image of this curve and close enough to $z$ : there exists a $t>0$ such that $x := \alpha(t) \in B(z,\frac{\delta_z}{2})\cap \Omega_0$ and the John condition yields that $\delta(x) \geqslant C^{-1}t$. Then, because $\alpha$ is a Lipschitz curve, it implies 
\[|x-z|=|\alpha(t)-\alpha(0)|\leqslant Ct \leqslant C^2 \delta(x).\]
We can assume, without loss of generality that $M \geqslant C^2$. Therefore, we found  a point $x\in  B(z, M\delta(x))\cap B(z,\frac{\delta_z}{2})\cap \Omega_0$. Thus this point $x$ is in $\Omega_0 \setminus \overline{\Omega}$, and, because $z \in \partial \Omega_0$,
\[\delta(x) = dist(x, \R^N \setminus \Omega_0) \leqslant |x-z| \leqslant \frac{\delta_z}{2}.\]
Hence, for all $y \in B(x,\frac{\delta(x)}{M})$, 
\[|y-z| \leqslant |y-x| +|x-z| \leqslant \frac{\delta(x)}{M} + \frac{\delta_z}{2} < \delta_z.\]
This yields that $y \notin \overline{\Omega}$, thus 
\[B(x,\frac{\delta(x)}{M}) \subset \R^N \setminus \overline{\Omega}.\]
Finally, this implies 
\[f_\sharp(x) = \frac{|B(x,\frac{\delta(x)}{M}) \setminus \overline{\Omega}|}{\delta(x)^N} = \frac{|B(x,\frac{\delta(x)}{M})|}{\delta(x)^N} = \frac{\omega_N}{M^N}, \]
and this completes the proof of the claim 
\[Nf(z) \geqslant f_\sharp(x) = \frac{\omega_N}{M^N}. \]
Then, we can write 
\begin{align*}
\Haus(\partial \Omega_0 \setminus \overline{\Omega}) &= \int_{\partial \Omega_0 \setminus \overline{\Omega}} 1 d\Haus \\
&\leqslant C \int_{\partial \Omega_0 \setminus \overline{\Omega}}Nf(z) d\Haus(z) \\
&\leqslant C \int_{\partial \Omega_0}Nf(z) d\Haus(z).
\end{align*}
And Lemma~\ref{lem 7.12 DS} yields 
\[\int_{\partial \Omega_0}Nf(z) d\Haus(z) \leqslant C \int_{\Omega_0} |\nabla f| = C \Haus(\partial \Omega \cap \Omega_0),\]
for some other constant (still denoted by $C$) only depending on $\Omega_0$. This achieves the proof of the first inequality 
\[\Haus(\partial \Omega_0 \setminus \overline{\Omega}\cap \mathcal{D}) \leqslant \Haus(\partial \Omega_0 \setminus \overline{\Omega})\leqslant C \Haus(\partial \Omega \cap \Omega_0 ) = C \Haus(\partial \Omega \cap \Omega_0 \cap \mathcal{D}).\]

\textit{Step 2. } For the second inequality the proof is very similar. We consider $f = \chi_{\Omega \setminus \Omega_0}$ and we reason with $\mathcal{D} \setminus \overline{\Omega_0}$ instead of $\Omega_0$, which is also a John domain with regular boundary since $\Omega_0$ is supposed to be a bi-John domain with regular boundary. For this case, we will suppose that the center of $\mathcal{D} \setminus \overline{\Omega_0}$, denoted $z_1$ is in $(\mathcal{D}\setminus \overline{E_0})\setminus \overline{\mathcal{C}}$ and that there exists a radius $r_1 \geqslant \Tilde{C}^{-1}$ with $\Tilde{C}$ a constant, such that $B_1 := B(z_1,r_1) \subset (\mathcal{D}\setminus \overline{E_0})\setminus \overline{\mathcal{C}} \subset \mathcal{D} \setminus \overline{\Omega_0}$. 
 In this case, for $x \in \mathcal{D}\setminus \overline{\Omega_0}$,  it yields $m_1f = 0$, because the boundary conditions impose $B_1 \subset (\mathcal{D}\setminus \overline{E_0})\setminus \overline{\mathcal{C}} \subset \mathcal{D} \setminus \overline{\Omega_0}$ and $B_1 \subset (\mathcal{D}\setminus \overline{E_0})\setminus \overline{\mathcal{C}} \subset \mathcal{D}\setminus \overline{\Omega}$. Thus, in this case, since by definition of $\delta(x):= dist(x,\overline{\Omega_0})$, $B(x,\frac{\delta(x)}{M}) \subset \R^N \setminus \overline{\Omega_0}$, and it follows 
\[f_{\sharp}(x) =\frac{1}{\delta(x)^N} \left|B\left(x,\frac{\delta(x)}{M}\right) \cap \Omega \right|.\]
Then, there exists a constant $C>0$ such that for all $z \in \partial \Omega_0 \cap \Omega$, $Nf(z) \geqslant C^{-1}$. The proof of this claim is similar to \textit{step 1}. And we can write 
\begin{align*}
\Haus(\partial \Omega_0 \cap \Omega) &= \int_{\partial \Omega_0 \cap \Omega} 1 d\Haus \\
&\leqslant C \int_{\partial \Omega_0 \cap \Omega}Nf(z) d\Haus(z) \\
&\leqslant C \int_{\partial\Omega_0}Nf(z) d\Haus(z).
\end{align*}
Moreover, Lemma~\ref{lem 7.12 DS}, applied with $W = \mathcal{D} \setminus \Omega_0$ yields 
\[\int_{\partial\Omega_0}Nf(z) d\Haus(z) \leqslant C \int_{\mathcal{D} \setminus \Omega_0} |\nabla f| = C \Haus(\partial \Omega \setminus \overline{\Omega_0} \cap \mathcal{D}),\]
for some other constant (still denoted by $C$) only depending on $\Omega_0$. This achieves the proof of the second inequality 
\[\Haus(\partial \Omega_0 \cap \Omega\cap \mathcal{D}) \leqslant \Haus(\partial \Omega_0 \cap \Omega)\leqslant C \Haus(\partial \Omega \setminus \overline{\Omega_0} \cap \mathcal{D}).\]

\end{proof}

\section{Example of a Lipschitz graph in a cylinder} 
\label{section cylindre}
 In this Section, we present an example of Plateau's problem where the boundary condition satisfies all Hypotheses~\ref{H2},\ref{H3},\ref{H1} required for the study of the regularity of P-quasi-minimizers. 

\subsection{Definition of the problem}

As convex set, we consider a cylinder $\mathcal{C}$ with convex basis $\mathcal{B}$, defined by \[\mathcal{C} = \mathcal{B} \times [-h,h] \subset \R^{N-1}\times \R.\]

Following the standard introduction of Plateau's problem, we introduce a prescribed curve $\Gamma$ as a Lipschitz graph contained in the boundary of the cylinder $ \partial \mathcal{C}$ (up to considering a higher cylinder, i.e. a larger $h>0$). Since we have assumed that $\Gamma$ is a Lipschitz graph, we denote by $\gamma : \partial \mathcal{B} \rightarrow [-h,h]$ the Lipschitz function, with Lipschitz constant $Lip(\gamma)$, such that
\[\Gamma = \{(\omega,\gamma(\omega)), \omega \in \partial \mathcal{B}\}.\]

\begin{figure}[ht]
    \centering
    \includegraphics[scale = 0.4]{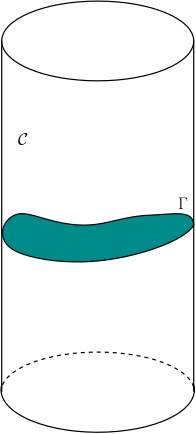}
    \caption{Example of surface in the cylinder $\mathcal{C}$ with boundary $\Gamma$.}
    \label{fig:cylinder}
\end{figure}

\begin{figure}[ht]
  \centering
    \subfloat[\centering Dilated cylinder]{\includegraphics[scale = 0.35]{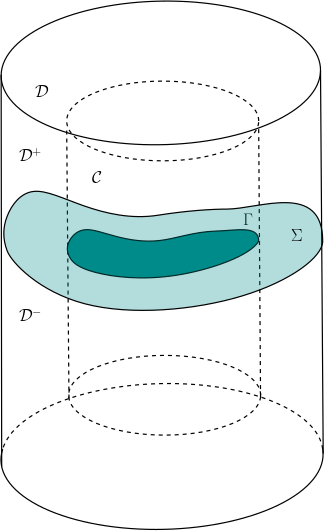} \label{fig:2 cylindes emboités}} \quad \quad \quad \quad \quad
    \subfloat[\centering Section of the cylinders]{\includegraphics[scale = 0.17]{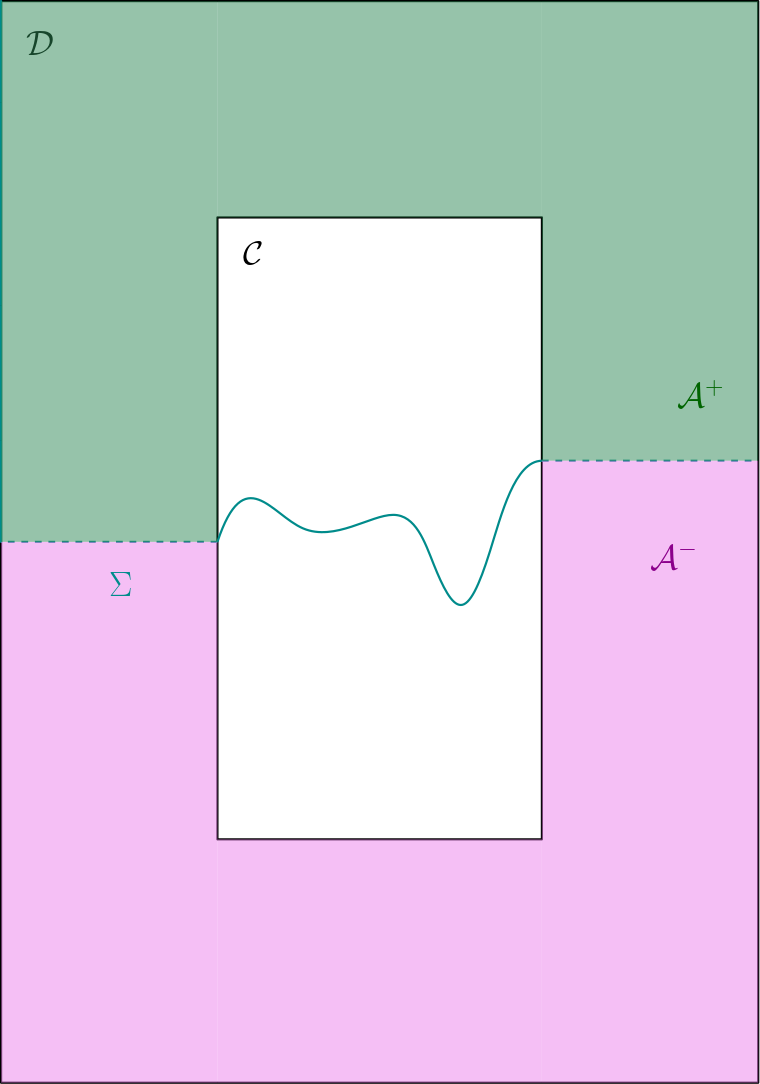} \label{fig:3 cylindes emboités}} 
  \caption{Definition of the boundary condition in a cylinder.}
  \label{fig: cylindes emboités}
\end{figure}

As bounded domain we use $\mathcal{D}$ the cylinder $\mathcal{C}$ which has been dilated by a factor $\lambda > 1$. In other words $\mathcal{D}=\lambda \mathcal{C}$, as shown in Figure~\ref{fig: cylindes emboités}, so that $\overline{\mathcal{C}} \subset \mathcal{D}$.

In our context, we need to define the boundary condition $E_0$ corresponding to the topological constraint to span the prescribed curve $\Gamma$. 
\begin{definition}
    We also dilate $\Gamma$ by the same factor $\lambda$, in $\mathcal{D}\setminus \overline{\mathcal{C}}$, and we define 
    \begin{equation}
    \label{def sigma}
    \Sigma := \{ \ (r\omega,\gamma(\omega)), r\in ]1,\lambda[ \text{ and } \omega \in \partial \mathcal{B} \ \}.
    \end{equation}
\end{definition}

Notice that the surface $\Sigma$ is a Lipschitz graph in $\mathcal{D}\setminus \overline{\mathcal{C}}$ that separates this region in two parts: one above the surface, denoted by 
\begin{align*}
    &\mathcal{A}^{+} := \{\ (r\omega,z), r \in ]1,\lambda[, \omega \in \partial \mathcal{B} \text{ and } z \in ]\gamma(\omega),h[ \ \},
\end{align*}
and one below, denoted by 
\[\mathcal{A}^{-} := \{\ (r\omega,z), r \in ]1,\lambda[, \omega \in \partial \mathcal{B} \text{ and } z \in ]-h,\gamma(\omega)[ \ \}.\]
We will minimize among sets, which include $\mathcal{A}^{+}$ (the part above the surface $\Sigma$) and do not intersect $\mathcal{A}^{-}$ (the part below $\Sigma$), namely whose complementary sets include $\mathcal{A}^{-}$. Thus, the boundary of such sets in the region $\mathcal{D}\setminus \overline{\mathcal{C}}$ coincides with the Lipschitz graph $\Sigma$, which translates the topological constraint of Plateau's problem. Thus, we define the boundary condition $E_0 = \mathcal{A}^+$. Notice that the boundary of $\mathcal{A}^+$ is regular in the sense that $\partial \mathcal{A}^+ =  \partial^* \mathcal{A}^+ $ and $\Sigma$, defined in Definition~\eqref{def sigma}, coincides with $\partial^* \mathcal{A}^+\cap \mathcal{D}$.

\subsection{Proof that the assumptions are satisfied}

As $\mathcal{A}^\pm$ are Lipschitz domains it is already known that they are 
domains of isoperimetry. Therefore, we will justify only \ref{H2} and \ref{H3}. 

\subsubsection{Ahlfors regularity of $\Sigma$}

We start by justifying that Hypothese~\ref{H2} that $\Sigma$ is Ahlfors regular. This comes from the classical result that a Lipschitz graph is Ahlfors regular.

\begin{prop}[Exterior Ahlfors regularity]
\label{th : ARext}
There exist some constants $C_{1},C_{2} > 0$, depending only on the dimension $N$ and the Lipschitz constant of the graph $\Sigma$, such that for all $x \in \Sigma$ and $r>0$ such that $B(x,r)\subset \mathcal{D}$,
    \begin{equation}
        C_{1} \leqslant \frac{\Haus(\Sigma \cap B(x,r))}{r^{N-1}} \leqslant C_{2}.
    \end{equation}
\end{prop}

\begin{proof}
Since $\Sigma$ is a Lipschitz graph, (we denote the associated Lipschitz function by $g$) there exists a set $A(x,r) \subset \R^{N-1}$ such that 
\[\Sigma\cap B(x,r)= \{(y_1,\cdots , y_{N-1}, g(y_1,\cdots , y_{N-1})) \ | \ (y_1,\cdots , y_{N-1}) \in A(x,r) \}.\]
Thus, the area formula (see for instance \cite[Theorem 2.71]{ambrosio2000oxford}) yields 
\[\Haus(\Sigma\cap B(x,r)) = \int_{A(x,r)}\sqrt{1+|\nabla g|^2}.\]
Then, $1 \leqslant \sqrt{1+|\nabla g|^2} \leqslant \sqrt{1+Lip(g)^2}$ implies 
\[|A(x,r)| \leqslant \Haus(\Sigma \cap B(x,r)) \leqslant \sqrt{1+Lip(g)^2} |A(x,r)|.\]
Notice that $ A(x,r) \subset B(\Tilde{x},r) $, where $\Tilde{x} = (x_1, \cdots, x_{N-1})$. Hence, 
\[ |A(x,r)| \leqslant \omega_{N-1}r^{N-1}.\]

Moreover, if $B(x,r) \subset \mathcal{D}\setminus \overline{\mathcal{C}}$ then 
$B\left(\Tilde{x}, \frac{r}{ \sqrt{1+Lip(g)^2}} \right) \subset A(x,r) $, which yields 
\[\omega_{N-1}\frac{1}{(\sqrt{1+Lip(g)^2})^{N-1}}r^{N-1} \leqslant |A(x,r)|.\]
Otherwise, since $\Sigma$ is a radial extension of $\Gamma$, there exists $y \in \Sigma \cap B(x,r/2)$ such that $B(y,r/2) \subset B(x,r)\setminus \overline{\mathcal{C}}$, in which case $B\left(\Tilde{y}, \frac{r}{2 \sqrt{1+Lip(g)^2}} \right) \subset A(y,r/2) \subset A(x,r) $. Thus,
\[\omega_{N-1}\frac{1}{(2\sqrt{1+Lip(g)^2})^{N-1}}r^{N-1} \leqslant |A(x,r)|.\]
This concludes the proof, with \[C_1 = \omega_{N-1}\frac{1}{(2\sqrt{1+Lip(g)^2})^{N-1}} \text{ and }C_2 = \omega_{N-1}\sqrt{1+Lip(g)^2}.\] 
\end{proof}

\subsubsection{Condition B'}
 In this section, we justify why Hypothesis~\ref{H3}, Condition B', is satisfied for the topological boundary of a competitor in $\mathcal{D}\setminus \overline{\mathcal{C}}$. Since, by definition, the boundary in this region is a Lipschitz graph, Proposition~\ref{prop : rect ext} only states that a Lipschitz graph satisfies Condition B'. Although this result is known, we recall the argument for the reader convinience. Notice also that in this region any competitor is open.
 
\begin{prop}
\label{prop : rect ext}
There exists a constant $C>1$ such that for all $x\in \Sigma $ and $r>0$ such that $B(x,r) \subset \mathcal{D}$ and two balls $B_1\subset (E_0 \cap B(x,r)) \setminus \overline{\mathcal{C}}$ and $B_2 \subset (B(x,r) \setminus \overline{E_0}) \setminus \overline{\mathcal{C}}$, with radii larger than $C^{-1}r$.
\end{prop}

\begin{proof}
Let $x = (x_1,...,x_{N-1}, g(x_1,...,x_{N-1})) \in \Sigma$ and $r>0$ be such that $B(x,r) \subset \mathcal{D}$. First we assume that $  B(x,r) \subset \mathcal{D}\setminus \overline{\mathcal{C}}$. 

We introduce the definition of two balls 
\[B_{1} = B(y^+, \rho), \ B_{2} = B(y^-, \rho),\]
where, 
\[y^{\pm} := (x_1,...,x_{N-1},g(x_1,...,x_{N-1}) \pm \alpha),\] 
and \[ \alpha := r\left(1-\frac{1}{2+Lip(g)}\right)>0 \ \text{ and } \ \rho := r\left(\frac{1}{2+Lip(g)}\right)>0 .\]
Naturally, both balls, designated as $B_1$ and $B_2$, are contained within the region $B(x,r) = B(x,r)\setminus \overline{\mathcal{C}}$. This is based on the specified values of $y_{\pm}$, $\alpha$, and the radius of the balls $\rho$. It thus remains to justify why $B_1$ is a subset of $E_0$ and why $B_2$ is a subset of $\mathcal{D}\setminus \overline{E_0}$. Since $B(x,r) \subset \mathcal{D}\setminus \overline{\mathcal{C}}$, it is equivalent to showing that they do not intersect $\Sigma$. Notice that $\Sigma$ separates the two regions $E_0 = \mathcal{A}^+$ and $\mathcal{D}\setminus \overline{E_0} = \mathcal{A}^-$. Moreover, because  $y^{+} \in \mathcal{A}^{+}$ and $y^{-}\in \mathcal{A}^{-} $, if the balls do not intersect with the hypersurface $\Sigma$, it can be concluded that they are, in fact, located within respectively $E_0$ and $\mathcal{D}\setminus \overline{E_0}$. So we will prove by contradiction that the balls do not cross $\Sigma$. Let us assume, by contradiction, that there exists a point $z$ belonging to $B_1 \cap \Sigma$. Then one can write $z = (z_1,..., z_{N-1}, g(z_1,...,z_{N-1}))$. And we have 
\begin{align*}
|g(z_1,...,z_{N-1}) - g(x_1,...,x_{N-1})| &= |g(z_1,...,z_{N-1}) - g(x_1,...,x_{N-1}) + \alpha - \alpha| \\
& \geqslant \alpha - | g(x_1,...,x_{N-1}) + \alpha - g(z_1,...,z_{N-1})| \\
& = \alpha -|y^{+}_{N}-z_N|\geqslant \alpha - |y^{+}-z|.
\end{align*}
But since we have assume that $z \in B_1$, it yields, by definition of $\alpha$ and $\rho$, 
\[|g(z_1,...,z_{N-1}) - g(x_1,...,x_{N-1})| > \alpha - \rho = Lip(g)\rho > Lip(g)|z'-y'| = Lip(g)|z'-x'|,\]
where $z' = (z_1,...,z_{N-1})$ and $y' = (y_1,...,y_{N-1}) = (x_1,...,x_{N-1}) = x'$. This is inconsistent with the assumption that $g$ is Lipschitz with constant $Lip(g)$. It therefore follows that $B_1$ does not intersect the surface $\Sigma$. The same line of reasoning shows that $B_2$ is also excluded from the surface $\Sigma$.

Secondly, if $B(x,r) \not \subset \mathcal{D}\setminus \overline{\mathcal{C}}$ then, because $\Sigma$ is a radial extension there exists $y \in \Sigma \cap B(x,r/2)$ such that $B(y,r/2) \subset B(x,r)\setminus \overline{\mathcal{C}}$. And appling the same computation as before in $B(y,r/2)$ leads to 
\[B_1 \subset B(y,r/2) \cap E_0 \subset (B(x,r) \cap E_0) \setminus \overline{\mathcal{C}} \] and 
\[B_2 \subset B(y,r/2) \setminus \overline{E_0} \subset (B(x,r) \setminus \overline{E_0}) \setminus \overline{\mathcal{C}},\]
 which concludes the proof.
\end{proof}

\section{Generalization to bi-Lipschitz image of a convex set}
\label{section BiLip}

The objective of this section is to justify how the optimal regularity results, which appears to be applicable solely for Plateau-quasi-minimizer within a convex set, can be extended to include Plateau's solutions in disparate geometric scenarios. More precisely we consider bi-Lipschitz images of the convex set $\mathcal{C}$. We begin by recalling the definition of a bi-Lipschitz image of a set, as follows:
\begin{definition}
$E$ is a bi-Lipschitz image of $F$ if there exit some Lipschitz and bijective application $\phi : F \longrightarrow E$ such that $E = \phi(F)$ and $\phi^{-1}$ is also Lipschitz.
\end{definition}

\begin{ex}
A cylinder with a curved axis and squared basis, denoted by $E$, is a bi-Lipschitz image of the rectangular box, denoted by $\mathcal{C}$(see Figure~\ref{fig bilip}). 
\begin{align*}
\phi :& \mathcal{C} \longrightarrow \R^{3}\\
		&(x,y,z) \longmapsto ((x-h)\cos(\frac{z \pi}{4h}),y,(x-h)\sin(\frac{z \pi}{4h}))
\end{align*} is indeed bi-Lipschitz.

\begin{figure}[ht]

\centering
\includegraphics[scale=0.2]{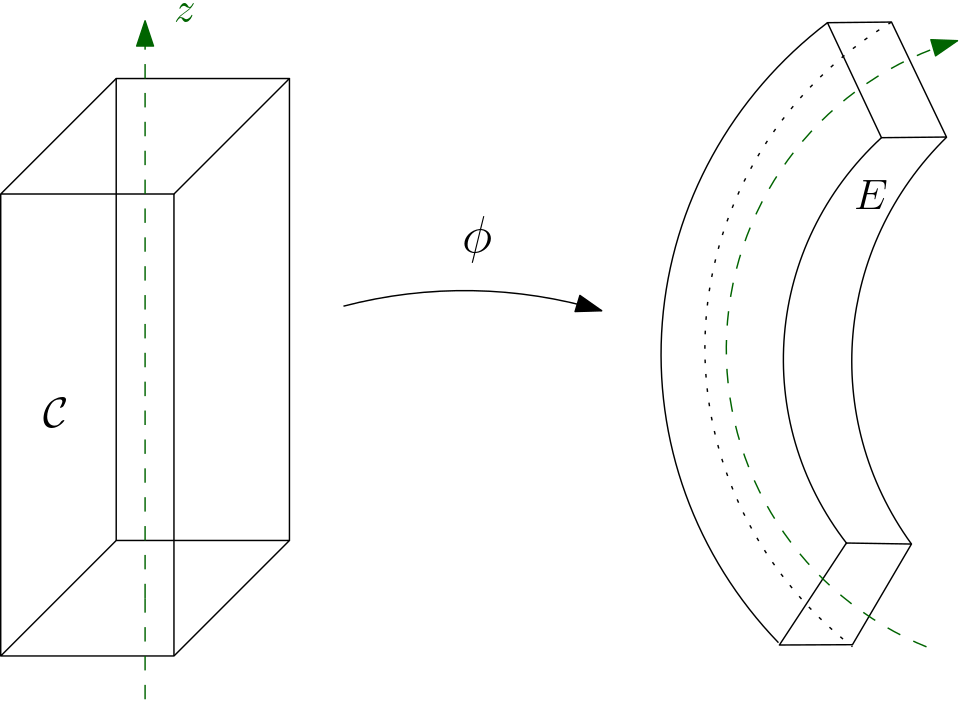}
\caption{Example of a bi-Lipschitz image of the squared based cylinder}
\label{fig bilip}
\end{figure}
\end{ex}

\subsection{Ahlfors regularity}

\begin{lemme}
\label{lem : 5.7}
Let $F$ and $G$ be some sets and $\phi: G \to F$ be a bi-Lipschitz application. We introduce the set $A \subset  \subset G$ of finite perimeter. Then we have
\[\partial ^{*} \phi(A) = \phi(\partial ^{*}A).\]
\end{lemme}

\begin{proof}

\cite[Theorem 1]{buczolich1992density} states that a bi-Lipschitz function maps the density points of $A$ into the density points of the bi-Lipschitz image $\phi(A)$. This directly yields the Lemma.
\end{proof}

 Now, fix a bi-Lipschitz image $\mathcal{F}$ of $\mathcal{D}$. Namely there exists a bi-Lipschitz function $\phi : \mathcal{D} \to \mathcal{F}$. And we introduce some notations : 
\[  \mathcal{F} := \phi(\mathcal{D}) \text{ and } \ \mathcal{E} := \phi(\mathcal{C}).\]

\begin{definition}
Consider Plateau's problem in $\mathcal{F}$, defined as follows
\begin{equation}
\label{pb limite image}
\inf \{P(\Omega,\mathcal{F}), \Omega \text{ such that } \chi_{\Omega} \in BV(\mathcal{F}), \Omega \setminus \phi(\overline{\mathcal{C}})= \phi(E_0) \setminus \phi(\overline{\mathcal{C}})\}.
\end{equation}
\end{definition}

\begin{rmq}
Notice first that for this Plateau problem, the boundary constraint is $\phi(E_0)$. Then the existence of solutions for this problem can easily be justified through the direct method of Calculus of Variations, as it is done in Proposition~\ref{existence}.
\end{rmq}

\begin{definition}[Plateau-quasi-minimizer of \eqref{pb limite image}]
A competitor for \eqref{pb limite image} $\Omega_{0}$ (\emph{i.e.}, $\chi_{\Omega_0} \in BV(\mathcal{F})$, $\Omega_0 \setminus \phi(\overline{\mathcal{C}})= \phi(E_0) \setminus \phi(\overline{\mathcal{C}})$) is called a Plateau-quasi-minimizer of \eqref{pb limite image} (or P-quasi-minimizer of \eqref{pb limite image}, for short), if there exists a constant $Q\geqslant 1$, such that for all $\Omega$ competitor for \eqref{pb limite image} we have
\[\Haus((\partial^*\Omega_0 \setminus \partial^*\Omega)\cap \mathcal{F})\leqslant Q \Haus((\partial^*\Omega \setminus \partial^*\Omega_0)\cap \mathcal{F}).\]
\end{definition}

\begin{lemme}
If $\Omega_{0}$ is a solution of \eqref{pb limite image} then $\phi^{-1}(\Omega_{0})$ is a Plateau-quasi-minimizer in the convex set $\mathcal{C}$ (see Definition~\ref{def quasi min}).
\end{lemme}

\begin{proof}
Let $\Omega_{0}$ be a solution of \eqref{pb limite image}. Let $\Omega$ be a competitor for Plateau's problem \eqref{def : pblimit}. First, let us check that $\phi(\Omega)$ is a competitor for \eqref{pb limite image}. By bijectivity and continuity of $\phi$ and because $\Omega$ is a competitor in $\mathcal{D}$, the condition $\phi(\Omega) \setminus \phi(\overline{\mathcal{C}})= \phi(E_0) \setminus \phi(\overline{\mathcal{C}})$ is satisfied. Besides, since $\phi^{-1}$ is Lipschitz and $\chi_{\phi(\Omega)}=\chi_{\Omega} \circ \phi^{-1}$ from \cite[Proposition 3.2]{ambrosio2000oxford} we get that the last condition $\chi_{\phi(\Omega)} \in BV(\mathcal{F})$. Hence, $\phi(\Omega)$ is a competitor for \eqref{pb limite image} and therefore
\begin{equation}
\label{equation de min}
\Haus((\partial^*(\phi(\Omega))\setminus \partial^* \Omega_0)\cap \mathcal{F}) \geqslant \Haus((\partial^*\Omega_0\setminus \partial^*( \phi(\Omega)))\cap \mathcal{F}).
\end{equation}
But, with Lemma~\ref{lem : 5.7}, and because $\phi$ is bi-Lipschitz, we have 
\begin{align*}
\Haus((\partial^*(\phi^{-1}(\Omega_{0}))\setminus \partial^* \Omega)\cap\mathcal{D})& = \Haus(\phi^{-1}((\partial^* \Omega_0 \setminus \partial^*(\phi( \Omega))) \cap \mathcal{F})) \\
&\leqslant Lip(\phi^{-1}) \Haus((\partial^* \Omega_0 \setminus \partial^*(\phi( \Omega))) \cap \mathcal{F}) \\
&\leqslant Lip(\phi^{-1}) \Haus((\partial^*(\phi(\Omega))\setminus \partial^* \Omega_0)\cap \mathcal{F}) \\
&\leqslant Lip(\phi^{-1})Lip(\phi) \Haus((\partial^* \Omega\setminus \partial^*(\phi^{-1}(\Omega_{0})))\cap\mathcal{D}).
\end{align*}
Thus, $\phi^{-1}(\Omega_{0})$ is a P-quasi-minimizer, with constant $Lip(\phi^{-1})Lip(\phi)$.
\end{proof}

\begin{rmq}
Conversely, by replacing $\phi$ by $\phi^{-1}$ in the previous proof, we get that if $\Omega_{0}$ is a minimizer of Plateau's problem \eqref{def : pblimit}, then $\phi(\Omega_{0})$ is a Plateau-quasi-minimizer of \eqref{pb limite image}.
\end{rmq}

\begin{lemme}
\label{prop : 4.8}
$\Omega_{0}$ is a Plateau-quasi-minimizer of \eqref{pb limite image} if and only if $\phi^{-1}(\Omega_{0})$ is a Plateau-quasi-minimizer of \eqref{def : pblimit}.
\end{lemme}

\begin{proof}
The proof of this Lemma is identical to the previous one. We just need to add the constant $Q$ in \eqref{equation de min}. Thus the constant $Lip(\phi)Lip(\phi^{-1})$ is replaced by $Q Lip(\phi)Lip(\phi^{-1})$.

\end{proof}

\begin{lemme}
\label{prop : 4.9}
Let $F$ and $G$ be some sets, for any bi-Lipschitz application $\phi: G \to F$, and $\Omega$ be a subset of $G$. Then, $\partial^*\Omega$ is Ahlfors regular, in the sense of Definition~\ref{def AR} if and only if $\partial^*\phi(\Omega)$ is Ahlfors regular.

\end{lemme}

\begin{proof}
Let $\phi$ be a bi-Lipschitz application from $G$ to $F$ and let $\Omega \subset G$.
Suppose that $\partial^*\Omega$ is Ahlfors regular. Let $x\in \partial^{*}\phi(\Omega)$ and $r>0$. From Lemma~\ref{lem : 5.7} we get that $x\in \phi(\partial^{*}\Omega)$, such that $\phi^{-1}(x) \in \partial^*\Omega$. Because $B(x,r) \subset \phi(B(\phi^{-1}(x), Lip(\phi^{-1})r)$, this yields
\begin{align*}
P(\phi(\Omega),B(x,r)) &\leqslant \Haus(\partial^{*}\phi(\Omega)\cap \phi(B(\phi^{-1}(x), Lip(\phi^{-1})r)))\\
 &= \Haus(\phi(\partial^{*}\Omega \cap B(\phi^{-1}(x), Lip(\phi^{-1})r))) \\ &\leqslant Lip(\phi)P(\Omega,B(\phi^{-1}(x), Lip(\phi^{-1})r))\\
&\leqslant Lip(\phi)Lip(\phi^{-1})^{N-1}C_{2}r^{N-1}.
\end{align*}
Using the same argument with $\phi^{-1}$, since $\phi(B(\phi^{-1}(x),\frac{r}{Lip(\phi)}))$, we get
\begin{align*}
\frac{C_{1}}{Lip( \phi)^{N-1}}r^{N-1} &\leqslant P(\Omega,B(\phi^{-1}(x),\frac{r}{Lip(\phi)})) \\
&= P(\phi^{-1}(\phi(\Omega)),\phi^{-1}(\phi(B(\phi^{-1}(x),\frac{r}{Lip(\phi)}))) \\
&\leqslant Lip(\phi^{-1})P(\phi(\Omega),\phi(B(\phi^{-1}(x),\frac{r}{Lip(\phi)})))\\
&\leqslant Lip(\phi^{-1})P(\phi(\Omega),B(x,r)).
\end{align*}
Hence, $\partial^*\phi(\Omega)$ is Ahlfors regular, with constants depending only on the dimension $N$, the Ahlfors constants of $\partial^*\Omega$ and the Lipschitz constants $Lip(\phi)$ and $Lip(\phi^{-1})$.  

By employing the same line of reasoning with respect to $\phi^{-1}$, we arrive at the converse implication.
\end{proof}

\begin{prop}
\label{th : 4.10}

 Let $\Omega_{0}$ be a Plateau-quasi-minimizer of \eqref{pb limite image}such that $\spt \mu_{\Omega_0} = \partial \Omega_0$. Then $ \partial^{*}\Omega_{0}$ is Ahlfors regular.
\end{prop}

\begin{proof}
Let $\Omega_{0}$ be a P-quasi-minimizer of \eqref{pb limite image}. Lemma~\ref{prop : 4.8} implies that $\phi^{-1}(\Omega_{0})$ is a P-quasi-minimizer of \eqref{def : pblimit}. Thus thanks to Proposition~\ref{th : AhlReg} we know that $\partial^*\phi^{-1}(\Omega_0)$ is Ahlfors regular. And then Lemma~\ref{prop : 4.9}, yields the Ahlfors regularity for the P-quasi-optimal surface in $E$ $\partial^*\Omega_{0}$.
\end{proof}

\subsection{Uniform rectifiability}

\begin{lemme}
\label{prop : bplg}
Let $F$ and $G$ be some sets, $\phi$ a bi-Lipschitz from $F$ to $G$, and $\Omega$ a subset of $F$. Then, $\Omega$ has BPLG if and only if $\phi(\Omega)$ has BPLG. 

\end{lemme}

\begin{proof}
The proof of this result is very similar to Lemma~\ref{prop : 4.9}. The only additional argument here is that the Lipschitz image of a Lipschitz graph is also a Lipschitz graph.
\end{proof}

\begin{prop}
\label{th : quasimin bplg}

 Let $\Omega_{0}$ be a Plateau-quasi-minimizer of \eqref{pb limite image}. Then, up to considering an equivalent set, $\partial \Omega_0$ has BPLG, so in particular is uniformly rectifiable.
\end{prop}
\begin{proof}
From Lemma~\ref{prop : 4.8} we know that $\phi^{-1}(\Omega_{0})$ is a P-quasi-minimizer. Thus, Proposition~\ref{th : unif rect} leads to $\phi^{-1}(\Omega_0)$ having BPLG (up to considering an equivalent set). And then Lemma~\ref{prop : bplg} yields that the essential boundary of the P-quasi-minimizer in $E$ has BPLG, so in particular is uniformly rectifiable.
\end{proof}

\subsection{Charectization by bi-John domain}

\begin{lemme}
\label{prop : john image}
Let $F$ and $G$ be open sets, $\phi$ a bi-Lipschitz application from $F$ to $G$, and $\Omega$ a subset of $F$. Then, $\Omega$ is a bi-John domain in $F$ if and only if $\phi(\Omega)$ is a bi-John domain in $G$. 

\end{lemme}

\begin{proof}
This is clear: we only need to compose the path by $\phi$, which is bi-Lipschitz. Hence, we get the desired equivalence. 
\end{proof}

\begin{thm}

$\Omega_{0}$ is a Plateau-quasi-minimizer of \eqref{pb limite image} if and only if there exist an open equivalent set, still denoted by $\Omega_0$ such that it is a bi-John domain with regular boundary, in the sense of Definition \ref{def : regbdd}, and a competitor for \eqref{pb limite image}.
\end{thm}
\begin{proof}
The Theorem follows directly from Lemma~\ref{prop : 4.8}, Lemma~\ref{prop : john image} and Theorem~\ref{th john dom}. Notice that Hypotheses~\ref{H2},\ref{H3},\ref{H1} are stable by bi-Lipschitz image.
\end{proof}

\bibliographystyle{alpha}
\bibliography{biblioarticle}

\end{document}